\documentclass[12pt]{article}
\usepackage{amssymb}
\usepackage{amsfonts}
\usepackage{times}
\usepackage{mathptmx}
\usepackage{amsmath}
\usepackage[usenames]{color}
\usepackage{mathrsfs}
\usepackage{amsfonts}
\usepackage{amssymb,amsmath}
\usepackage{cite}
\usepackage{cases}
\usepackage{amsthm}
\def\cl{\centerline}

\def\a{\alpha}

\def\vs{\vspace*}
\def\W{\mathscr{W}}

\def\G{\mathscr{G}}

\def\Z{\mathbb{Z}}

\def\C{\mathbb{C}}

\def\Der{{\rm {Der}}}

\def\vs{\vspace*}

\numberwithin{equation}{section}
\newtheorem{theo}{Theorem}[section]
\newtheorem{defi}[theo]{Definition}
\newtheorem{coro}[theo]{Corollary}
\newtheorem{lemm}[theo]{Lemma}
\newtheorem{prop}[theo]{Proposition}

\newtheorem{ex}[theo]{Example}
\newtheorem{case}{Case}

\newtheorem{remark}[theo]{Remark}

\allowdisplaybreaks[4]
\begin{document}
\begin{center}
    {\bf\large Biderivations of Hom-Lie algebras and superalgebras}
    \footnote {Corresponding author: lmyuan@hit.edu.cn (Lamei Yuan)
    }
\end{center}

\cl{Lamei Yuan, Jiaxin Li}

\cl{\small School of Mathematics, Harbin Institute of Technology, Harbin
    150001, China}
\cl{\small E-mail: lmyuan@hit.edu.cn, 2683205187@qq.com}
\vs{6pt}

\vs{8pt}

{\small\footnotesize
\parskip .005 truein
\baselineskip 3pt \lineskip 3pt
\noindent{{\bf Abstract:} On Hom-Lie algebras and superalgebras, we introduce the notions of biderivations, linear commuting maps and $\alpha$-biderivations, and compute them for some typical Hom-Lie algebras and superalgebras, including $q$-deformed $W(2,2)$ algebra, $q$-deformed Witt algebra and superalgebra. 


\vs{10pt}

\noindent{\bf Key words:} Hom-Lie algebras, biderivations, linear commuting maps, $\alpha$-biderivations, $q$-deformed $W(2,2)$ algebra, $q$-deformed Witt algebra, $q$-deformed Witt superalgebra\\
\vs{5pt}
\noindent{\bf Mathematics Subject Classification (2000):} 17A30,
17A60, 17B05, 17B60, 17B70.

\parskip .001 truein\baselineskip 6pt \lineskip 6pt

\section{Introduction}
\vs{8pt}

The study of biderivations and commuting maps has its roots in associative ring theory\cite{B1,B2,P}. The development of the theory of such maps and their applications to different areas, in particular to
Lie theory, was surveyed in \cite{B3}. In the last decade, an enormous amount of attention has been put on the study
of biderivations and related maps on Lie algebras and superalgebras (see, for example, \cite{B4,LGZ,Chen2016,HWX, wang2013,Tang2018,Tang2017,Tang-Li2018,XWH,XS}).
\begin{defi}\label{d1}{\rm Let $L$ be a Lie algebra.
 A bilinear map $\varphi:L\times L\rightarrow L$ is called a biderivation if for every $x\in{L}$, the maps $y\mapsto\varphi(x,y)$ and  $y\mapsto \varphi(y,x)$ are derivations of ${L}$, i.e, for all $x,y,z\in L$, there hold
\begin{align}
 \varphi([x,y],z)&= [\varphi(x,z),y]+[x,\varphi(y,z)], \label{bi-1}\\
 \varphi(x,[y,z])&= [\varphi(x,y),z]+[y,\varphi(x,z)].\label{bi-2}
 \end{align}}
\end{defi}
For every element $\lambda$ in the ground field, the map $(x,y)\mapsto \lambda[x,y]$ is a biderivation on $L$.
Such maps are called {\it inner biderivations}.

Biderivations and commuting maps were studied by M. Bre$\rm\check{s}$ar in a series of papers \cite{B1,B2,B3,B4}. In particular, he and his coworker proposed a generalization of biderivations for Lie algebras involving modules in \cite{B4}.
\begin{defi}\label{d2}{\rm (\cite{B4}) Let $L$ be a Lie algebra and $M$ an $L$-module.
 A bilinear map $\delta:{L}\times {L}\rightarrow M$ is called a skew-symmetric biderivation if it satisfies $\delta (x,y)=-\delta (y,x)$ and
\begin{align}
\delta([x,y],z) &= x \cdot\delta(y,z) -y\cdot \delta(x,z), \label{d4}
 \end{align}
for all $x,y,z\in L$.}
\end{defi}
Notice that relation \eqref{d4} means that the map $x\mapsto\delta(x,z)$ is a derivation from ${L}$ to $M$ for every $z\in{L}$. And the skew-symmetry of $\delta$ implies that $\delta$ is also
a derivation with respect to the second component. If $M = L$ and $x\cdot y = [x,y]$ in Definition \ref{d2}, it reduces to the usual biderivations (see Definition \ref{d1}) satisfying the additional condition: skew-symmetry. Fortunately, in some Lie algebras, the inner biderivations (which are obviously skew-symmetric) are the only biderivations (see, for example,  \cite{Chen2016,wang2013,HWX}). Hence, a number of known results \cite{Chen2016,wang2013,HWX} can be easily deduced by the technique developed in \cite{B4}. However,
some Lie algebras, like the twisted Heisenberg-Virasoro algebra \cite{Tang-Li2018} and block Lie algebra \cite{LGZ}, have biderivations that are not skew-symmetric. It seems much safer to use Definition \ref{d1} to compute biderivations when no modules are involved.

The motivations to study Hom-Lie structures are related to physics and deformations of Lie
algebras, in particular Lie algebras of vector fields. The paradigmatic examples are $q$-deformations of
Witt and Virasoro algebras constructed in \cite{AS,CJ,H}.
The notion of Hom-Lie algebra was initially introduced in \cite{HLS} and then extended
to quasi-hom Lie and quasi-Lie algebras in \cite{LS1,LS2}. A slightly more general and commonly used definition of Hom-Lie algebras was given in
\cite{MS1}. The notion, constructions and properties of enveloping algebras and Chevalley-Eilenberg type homology of Hom-Lie algebras were studied in \cite{Yau1,Yau2}. Cohomology and deformations of Hom-Lie algebras were studied in \cite{AEM,MS2}.
The representation theory of Hom-Lie algebras that are multiplicative  (specifically, the twisting map is an algebra homomorphism) was developed in \cite{SH}. As a $q$-deformation of the centerless $W(2,2)$ Lie algebra \cite{ZD}, the $q$-deformed $W(2,2)$ algebra was constructed and studied in \cite{Yuan}, and it was proved to be a Hom-Lie algebra in \cite{YY}.

The Hom-Lie superalgebras were introduced in \cite{AM}, in which the $q$-deformed Witt superalgebra was constructed. The cohomology theory of Hom-Lie superalgebras was introduced and the second cohomology of the $q$-deformed Witt superalgebra was computed in \cite{AMS1,AMS2}. The notion of Hom-Lie color algebra was introduced and studied in \cite{Y1}.
It is a natural generalization of Hom-Lie algebras as well as a special
case of quasi-hom-Lie algebras.

Following \cite{B4}, Sun et al. introduced a Hom-generalization of biderivations in \cite{sun2020} and
proposed a definition of biderivations from a ``multiplicative" Hom-Lie algebra to its Hom-modules.
In their study, aside from restricting themselves to Hom-Lie algebras that are multiplicative, some results are only valid for biderivations of skew-symmetry. In this article, we aim to present a more general definition of biderivations without assuming skew-symmetry for arbitrary Hom-Lie algebras that are not restricted to be multiplicative. The natural choice comes from a Hom analogue of Definition \ref{d1}. This leads us to the following concept:
\begin{defi}\label{d3}{\rm Let $(\mathscr{G},[\cdot,\cdot],\alpha)$ be a Hom-Lie algebra, a bilinear map $\varphi:\mathscr{G}\times \mathscr{G}\rightarrow \mathscr{G}$ is called a biderivation of $\mathscr{G}$ if
        \begin{align}
            \varphi([x,y],\alpha(z))= & ~[\varphi(x,z),\alpha(y)]+[\alpha(x),\varphi(y,z)],\label{vf1} \\
            \varphi(\alpha(x),[y,z])= & ~[\varphi(x,y),\alpha(z)]+[\alpha(y),\varphi(x,z)],\label{vf2}
        \end{align}
        for all $x,y,z\in \mathscr{G}$.}
\end{defi}
We see that, if $\alpha={\rm id}$, then $\mathscr{G}$ is a Lie algebra and it reduces to Definition \ref{d1}.
Hence, we believe this is the only possible choice to give Hom-generalizations of biderivations if we do not need to take into account Hom-modules. We note a close connection between biderivations and linear commuting maps (see Definition \ref{d5}), and compute them for some typical examples of Hom-Lie algebras.  We will also consider the super case in this article.

The rest of this article is arranged as follows. In Section 2, we recall the basic definitions of Hom-Lie algebras and superalgebras, and write down main examples, including the $q$-deformed $W(2,2)$ algebra, the $q$-deformed Witt algebra and superalgebra. Moreover,
we introduce the notions of biderivations and linear commuting maps on Hom-Lie algebras and superalgebras, and characterize  their connections.

In Sections 3, we will compute biderivations on the $q$-deformed $W(2,2)$ algebra, the $q$-deformed Witt algebra and superalgebra, respectively. We observe that these Hom-algebras are $\Z$-graded, resulting in that the corresponding spaces of biderivations are $\Z$-graded. Hence, the inspection of biderivations can be performed by analyzing homogenous biderivations. Though the commutation looks like straightforward, it needs some techniques due to the $q$-bracket. Finally,
we characterize linear commuting maps on these algebras by applying the corresponding results on biderivations.

In Section 4,
we introduce the notions of $\alpha$-derivations and $\alpha$-biderivations for a Hom-Lie algebra $(\G,[\cdot,\cdot],\alpha)$, and we establish a close relation between them. The superanalogue is also studied. We will show that the $q$-deformed $W(2,2)$ algebra, $q$-deformed Witt algebra and superalgebra have no nonzero $\alpha$-biderivations or $\alpha$-super-biderivations.  We also give an example of Hom-Lie superalgebras with nonzero $\alpha$-super-derivations and $\alpha$-super-biderivations.

Throughout this paper, all the vector spaces are assumed to be over $\C$,
the field of complex numbers. For
a vector spaces $V$, we denote by ${\rm id}$ the identity map on $V$. Moreover, we denote by $\mathbb{Z}$ the ring of integers and
by $\mathbb{Z}_2 =\{\bar0,\bar1\}$ the cyclic group of order 2. For a superspace $V=V_{\bar0}\bigoplus V_{\bar1}$, we denote by $\mathcal{H}(V)$ the set of homogenous elements of $V$.




\section{Preliminaries}
\vs{8pt}

In this section, we summarize the definitions of Hom-Lie algebras and superalgebras. Also, we recall
the $q$-deformed W(2,2) algebra, the $q$-deformed Witt algebra and its supercase. Moreover, we introduce the notions of  biderivations and linear commuting maps of Hom-Lie algebras and their superanalogue, and we characterize connections between biderivations and linear commuting maps.


\subsection{Biderivations of Hom-Lie algebras}
\begin{defi}{\rm (\cite{MS1})
        A Hom-Lie algebra is a triple
        $(\mathscr{G},[\cdot,\cdot],\alpha)$ consisting of a vector space $\mathscr{G}$, a
        bilinear map $[\cdot,\cdot]:\mathscr{G}\times \mathscr{G}\rightarrow \mathscr{G}$ and a
        linear map $\alpha:\mathscr{G}\rightarrow \mathscr{G}$, such that
        \begin{align}
            [x,y]=                                                   & ~-[y,x], \ \ (\mbox{skew-symmetry}) \\
            {[\alpha(x),[y,z]]+[\alpha(y),[z,x]]+[\alpha(z),[x,y]]}= & ~0,\ \ \
            (\mbox{Hom-Jacobi identity})
        \end{align}
        for all $x,y,z\in \mathscr{G}$.}
\end{defi}
A Hom-Lie algebra $(\mathscr{G},[\cdot,\cdot],\alpha)$ is called a {\it multiplicative} Hom-Lie algebra if $\a$ is an algebraic homomorphism, i.e., for any $x,y\in\mathscr{G}$, we have $\a([x,y])=[\a(x),\a(y)]$.

A Hom-Lie algebra $(\mathscr{G},[\cdot,\cdot],\alpha)$ is said to be {\it $\mathbb{Z}$-graded} if the underlying vector space $\mathscr{G}$ is $\mathbb{Z}$-graded, i.e., $\mathscr{G}=\bigoplus_{n\in\mathbb{Z}}\mathscr{G}_n$, and if, furthermore, $[\mathscr{G}_m, \mathscr{G}_n]\subseteq \mathscr{G}_{m+n}$ for all $m,n\in\mathbb{Z}$ and $\alpha$ is an even map, i.e., $\alpha(\mathscr{G}_m)\subseteq \mathscr{G}_m$ for all $m\in\mathbb{Z}$.

\begin{ex}{\rm ({\bf The $q$-deformed $W(2,2)$ algebra})
        Let $q$ be a fixed nonzero complex number such that $q$ is not a root of unity. For $n\in\Z$, let $[n]$ denote
        the $q$-number
        \begin{eqnarray*}\label{qnumber}
            [n]=\frac{q^n-q^{-n}}{q-q^{-1}}.
        \end{eqnarray*}
        It is easy to see that the $q$-number satisfies the following properties:
        \begin{eqnarray*}
            [-n]=-[n],\ \ \, q^{n}[m]-q^{m}[n]=[m-n], \ \ \ q^{-n}[m]+q^{m}[n]=[m+n].\label{q}
        \end{eqnarray*}
        Let $\mathcal{W}$ be a vector space with basis $\{L_n,W_n|n\in\mathbb{Z}\}$ and satisfying
        \begin{eqnarray}\label{H-L1}
            [L_n,L_m]=[m-n]L_{m+n},\ \
            {[L_n,W_m]=[m-n]W_{m+n}},\ \forall~ m,n\in\Z,
        \end{eqnarray}
        and other brackets are obtained by skew-symmetry or are equal to 0.

        The linear map $\alpha$ on $\mathcal {W}$ is defined on the generators by
        \begin{eqnarray}\label{alpha1}
            \alpha(L_n)=(q^n+q^{-n})L_n,\ \  \alpha(W_n)=(q^n+q^{-n})W_n, \ \,\forall~ n\in\Z.
        \end{eqnarray}
        It was proved in \cite{YY} that the triple ($\mathcal {W}, [\cdot,\cdot], \alpha$)
        forms a Hom-Lie algebra, called the $q$-deformed $W(2,2)$ algebra and denoted by $\mathcal {W}^q$.  By defining ${\rm deg}(L_n)={\rm deg}(W_n)=n$, we obtain that $\W^q$ is a $\mathbb{Z}$-graded Hom-Lie algebra, i.e.,
        $$\W^q=\bigoplus_{n\in\mathbb{Z}}\W^q_n, \ \ \ \W^q_n={\rm span}_{\mathbb{C}}\{L_n,W_n\}.$$
        }
\end{ex}

\begin{ex}{\rm ({\bf The $q$-deformed Witt algebra}) Let $q\in\C\setminus\{0,1\}$ and $n\in\mathbb{Z}$. We set a second $q$-number as follows:
        \begin{eqnarray}\label{q-mumber}
            \{n\}=\frac{1-q^n}{1-q}.
        \end{eqnarray}
        The $q$-numbers have the following properties
        \begin{eqnarray*}
            \{n+1\}=1+q\{n\}=\{n\}+q^n,\ \,\{n+m\}=\{n\}+q^n\{m\},\ \,q^n\{-n\}=-\{n\}.
        \end{eqnarray*}
        Let $\mathcal {V}$ be a vector space with $\{L_n|n\in\mathbb{Z}\}$ as a basis. Define a bilinear operation on the generators by
        \begin{eqnarray*}
            [L_n,L_m]=(\{m\}-\{n\})L_{m+n}, \ \forall~  m,n\in\Z,
        \end{eqnarray*}
        and a linear map $\alpha$ by
        \begin{eqnarray*}
            \alpha(L_n)=(1+q^n)L_n, \forall~ n\in\Z.
        \end{eqnarray*}
        It is easy to check that $(\mathcal {V}, [\cdot,\cdot], \alpha)$
        is a Hom-Lie algebra, called the $q$-deformed Witt algebra and denoted by $\mathcal {V}^q$.
        It was constructed in \cite{HLS}
        as the first example of Hom-Lie algebras.}
\end{ex}
 It is easy to see that both $\W^q$ and $\mathcal {V}^q$ are not multiplicative. Recall that we have given the definition of biderivations on Hom-Lie algebras in Definition \ref{d3}. Let $(\mathscr{G},[\cdot,\cdot],\alpha)$ be a Hom-Lie algebra. An example of biderivations of $\mathscr{G}$ is the map
\begin{eqnarray}
    \varphi:\mathscr{G}\times \mathscr{G}\rightarrow \mathscr{G},~\varphi(x,y)=\lambda[x,y],~\forall~ x,y \in \mathscr{G}, \label{inner}
\end{eqnarray}
where $\lambda \in \C$. Such biderivations are called {\it inner biderivations}, and all other biderivations are called {\it outer biderivations}. We denote by ${\rm BDer}(\mathscr{G},\mathscr{G})$ and ${\rm BInn}\,(\mathscr{G},\mathscr{G})$ the vector space consisting of all biderivations and inner biderivations, respectively.
For convenience, we denote by $\varphi_{ad}$ the inner biderivation defined by \eqref{inner} with $\lambda=1$, which is called the {\it standard inner biderivation} on $\mathscr{G}$.

Let $(\mathscr{G},[\cdot,\cdot],\alpha)$ be a $\mathbb{Z}$-graded Hom-Lie algebra and $\varphi\in {\rm BDer}(\mathscr{G},\mathscr{G})$. For arbitrary $x,y\in \mathscr{G}$, $\varphi(x,y)$ will be a finite linear combination (depending on $x$ and $y$) of basis elements in $\mathscr{G}.$  We call a biderivation $\varphi$ {\it homogeneous of degree $d$} if there exists a $d\in\mathbb{Z}$ such that for all $i,j\in\Z$ and homogeneous elements $x_i,x_j\in\mathscr{G}$ of
${\rm deg}(x_i)=i$ and ${\rm deg}(x_j)=j$, we have

\begin{eqnarray*}
    \varphi(x_i,x_j)\in\mathscr{G}_n,\ \,\mbox{with}\ \, n=i+j+d.
\end{eqnarray*}
The corresponding subspace of homogeneous biderivations of degree $d$ is denoted by ${\rm BDer}_d(\mathscr{G},\mathscr{G})$.
Every biderivation $\varphi$ can be written as a formal infinite sum
\begin{eqnarray*}
    \varphi =\sum_{d \in \Z} \varphi_{(d)}, \ \ ~ \varphi_{(d)}\in {\rm BDer}_d(\mathscr{G},\mathscr{G}).
\end{eqnarray*}
Note that evaluated for a fixed pair of elements only a finite number of the summands will produce values different from zero.

\begin{defi}{\rm \label{d5} Let $(\mathscr{G},[\cdot,\cdot],\alpha)$ be a Hom-Lie algebra, and $f:\mathscr{G}\rightarrow \mathscr{G}$ a linear map. If $f\circ\alpha=\alpha\circ f$ and
        \begin{align}
            [f(x),x]=0,~\forall~ x \in \mathscr{G},\label{f(x)}
        \end{align}
        then $f$ is called a linear commuting map on $\mathscr{G}$.
    }\end{defi}

Obviously, condition \eqref{f(x)} is equivalent to the following
\begin{eqnarray}\label{L-M}
    [f(x),y]=[x,f(y)],~\forall~x,y\in \mathscr{G}.
\end{eqnarray}
Linear commuting maps and biderivations are intimately connected in the following way:
\begin{lemm}\label{lemma1}
    Let $f$ be a linear commuting map on a Hom-Lie algebra $(\mathscr{G},[\cdot,\cdot],\alpha)$. Define
    \begin{eqnarray}\label{L-M-1}
        \varphi(x,y)=[x,f(y)],~\forall~x,y \in \mathscr{G}.
    \end{eqnarray}
    Then $\varphi$ is a biderivation.
\end{lemm}
\begin{proof}
    By \eqref{L-M-1}, we have
    \[
        \varphi([x,y],\alpha(z))=\big[[x,y],f\big(\alpha(z)\big)\big].
    \]
    On the other hand, by Hom-Jacobi identity, we have
    \[
        [\varphi(x,z),\alpha(y)]+[\alpha(x),\varphi(y,z)]=\big[[x,f(z)],\alpha(y)\big]+\big[\alpha(x),[y,f(z)]\big]=\big[[x,y],\alpha\big(f(z)\big)\big].
    \]
    With $ f\circ\alpha=\alpha\circ f$, we get
    \[
    \varphi([x,y],\alpha(z))= [\varphi(x,z),\alpha(y)]+[\alpha(x),\varphi(y,z)].
    \]
    By \eqref{L-M} and \eqref{L-M-1}, we have
    \[
        \varphi(\alpha(x),[y,z])=\big[f\big(\alpha(x)\big),[y,z]\big].
    \]
    By Hom-Jacobi identity, we have
    \[
        [\varphi(x,y),\alpha(z)]+[\alpha(y),\varphi(x,z)]=\big[[f(x),y],\alpha(z)\big]+\big[\alpha(y),[f(x),z]\big]=\big[\alpha\big(f(x)\big),[y,z]\big].
    \]
    With $ f\circ\alpha=\alpha\circ f$, we obtain
    \[
        \varphi(\alpha(x),[y,z])=[\varphi(x,y),\alpha(z)]+[\alpha(y),\varphi(x,z)].
    \]
   This ends the proof.
\end{proof}

\subsection{Super-biderivations of Hom-Lie superalgebras}

Let $V$ be a linear superspace that is a $\mathbb{Z}_2$-graded linear space with a direct sum
$V = V_{\bar 0}\bigoplus V_{\bar1}$. The elements of $V_j$, $j=\{\bar0,\bar1\}$,  are said to be homogenous and of parity $j$. The parity of
a homogeneous element $x$ is denoted by $|x|$. Let $V$ and $W$ be two linear superspaces. A linear map from
$V$ and $W$ is said to be homogeneous of degree $s\in\mathbb{Z}_{2}$, if
\begin{align}
    f(V_j)\subset W_{j+s},\ \ ~ j=\{\bar0,\bar1\}.
\end{align}
If, in addition, $f$ is homogenous of degree zero, namely, $f(V_j)\subset W_{j},\ ~ j=\{\bar0,\bar1\}$,
then $f$ is said to be even. The degree of
a homogeneous linear map $f$ is denoted by $|f|$.

\begin{defi}{\rm (\cite{AM})
        A Hom-Lie superalgebra is a triple $(\mathscr{G},[\cdot,\cdot],\alpha)$ consisting of a superspace $\mathscr{G}$, an even bilinear map $[\cdot,\cdot]: \mathscr{G}\times \mathscr{G} \to \mathcal{G}$, and an even superspace homomorphism $\alpha:\mathscr{G} \to \mathcal{G}$ satisfying
        \begin{eqnarray}
            &&[x,y]=-(-1)^{|x||y|}[y,x],\\
            &&(-1)^{|x||z|}[\alpha(x),[y,z]]+(-1)^{|z||y|}[\alpha(z),[x,y]]+(-1)^{|y||x|}[\alpha(y),[z,x]]=0,\label{jacobi}
        \end{eqnarray}
        for all homogeneous elements $x, y, z$ in $\mathscr{G}$.}
\end{defi}

Equation \eqref{jacobi} is equivalent to the following
\begin{eqnarray}
    [\alpha(x),[y,z]]=[[x,y],\alpha(z)]+(-1)^{|x||y|}[\alpha(y),[x,z]],\label{jac}
\end{eqnarray}
for all homogeneous elements $x, y, z$ in $\mathscr{G}$.

\begin{defi}{\rm Let $(\mathscr{G},[\cdot,\cdot],\alpha)$ be a Hom-Lie superalgebra. A homogeneous linear map $D:\mathscr{G}\rightarrow \mathscr{G}$ is called a super-derivation of $\mathscr{G}$ if
        $$D([x,y])=[D(x),y]+(-1)^{|x||D|}[x,D(y)], \ \forall\ \, x\in \mathcal{H}(\mathscr{G}),\ y\in\mathscr{G}.$$
    }\end{defi}


A Hom-Lie superalgebra $(\mathscr{G},[\cdot,\cdot],\alpha)$ is called {\it multiplicative} if $\a$ is an endomorphism of Hom-Lie superalgebras.

\begin{ex}{\rm({\bf The $q$-deformed Witt superalgebra}) For $q\in\C\setminus\{0,1\}$ and $n\in\Z$, let $\{n\}$ be the $q$-number defined by \eqref{q-mumber}.
        A $q$-deformed Witt superalgebra, denoted by $\mathcal {L}^q$, can be presented as the $\Z_2$-graded vector space with $\{L_n|\,n\in \Z\}$ as a basis of the even homogeneous part and $\{G_n|\,n\in \Z\}$ as a basis of the odd homogeneous part. It is equipped with the commutator
        \begin{align}
            [L_n, L_m]= (\{m\}-\{n\})L_{n+m},   \
            [L_n, G_m]= (\{m+1\}-\{n\})G_{n+m},\ \forall \ m,n\in\Z.
        \end{align}
        The other brackets are obtained by supersymmetry or are equal to $0$.
        The even linear map $\alpha$ on $\mathcal {L}^q$ is defined on the generators by
        \begin{eqnarray}
            \alpha(L_n)=(1+q^n)L_n, \quad  \alpha(G_n)=(1+q^{n+1})G_n, \ \forall \ n\in\Z.
        \end{eqnarray}
        This algebra was constructed in \cite{AM} as an example of Hom-Lie superalgebras, and it is not multiplicative.}
\end{ex}

Now we introduce the notion of super-biderivations on Hom-Lie superalgebras.

\begin{defi}{\rm Let $(\mathscr{G},[\cdot,\cdot],\alpha)$ be a Hom-Lie superalgebra.
        A homogeneous bilinear map $\varphi:\mathscr{G}\times \mathscr{G}\rightarrow \mathscr{G}$ is said to be a super-biderivation of $\mathscr{G}$ if it satisfies

        \begin{align}
            \varphi \big([x,y],\alpha(z)\big)=~ & (-1)^{|y||z|}[\varphi(x,z),\alpha(y)]+(-1)^{|\varphi||x|}[\alpha(x),\varphi(y,z)],\label{svf1} \\
            \varphi\big(\alpha(x),[y,z]\big)=~  & [\varphi(x,y),\alpha(z)]+(-1)^{(|\varphi|+|x|)|y|}[\alpha(y),\varphi(x,z)],\label{svf2}
        \end{align}
        for all homogeneous elements $x, y, z$ in $\mathscr{G}$.}
\end{defi}

An example of super-biderivations of $\mathscr{G}$ is the map
\begin{eqnarray}
    \varphi:\mathscr{G}\times \mathscr{G}\rightarrow \mathscr{G},~\varphi(x,y)=\lambda[x,y],~\forall~ x,y \in \mathscr{G}, \label{s-inner}
\end{eqnarray}
where $\lambda \in \C$. Such super-biderivations are called {\it inner super-biderivations}, and all other super-biderivations are called {\it outer super-biderivations}. We still denote by ${\rm BDer}(\mathscr{G},\mathscr{G})$ and ${\rm BInn}\,(\mathscr{G},\mathscr{G})$ the vector space consisting of all super-biderivations and inner super-biderivations, respectively. In particular, we denote by $\varphi_{ad}$ the inner super-biderivation defined by \eqref{s-inner} with $\lambda=1$, which is called the {\it standard inner super-biderivation} on $\mathscr{G}$.

A super-biderivation $\varphi$ of a Hom-Lie superalgebra $(\mathscr{G},[\cdot,\cdot],\alpha)$ is called homogenous of degree $\gamma \in \Z_2$ if it is a
super-biderivation such that $\varphi(\mathscr{G}_\alpha,\mathscr{G}_\beta)\subset \mathscr{G}_{\alpha+\beta+\gamma}$ for any $\alpha,\beta\in \Z_2$. Denote by
${\rm BDer}_\gamma(\mathscr{G},\mathscr{G})$ the subspace of all homogenous super-biderivations of degree $\gamma$ of $\mathscr{G}$. Hence,
$${\rm BDer}(\mathscr{G},\mathscr{G}) = {\rm BDer}_{\bar{0}}(\mathscr{G},\mathscr{G})\oplus {\rm BDer}_{\bar{1}}(\mathscr{G},\mathscr{G}).$$

\begin{defi}{\rm Let $(\mathscr{G},[\cdot,\cdot],\alpha)$ be a Hom-Lie superalgebra. A homogeneous linear map $f:\mathscr{G}\rightarrow \mathscr{G}$ is called a linear commuting map if it satisfies
        \begin{align}
            f\circ\alpha=\alpha\circ f,\ \ ~
            [f(x),y]=-(-1)^{|x||y|}[f(y),x],\label{super-f(x)}
        \end{align}for all homogeneous elements $x, y$ in $\mathscr{G}$.
    }\end{defi}

We have the following result with a similar proof of Lemma \ref{lemma1}:
\begin{lemm}\label{lemma2}
    Let $f$ be a linear commuting map on a Hom-Lie superalgebra $(\mathscr{G},[\cdot,\cdot],\alpha)$. Define
    \begin{eqnarray*}
        \varphi(x,y)=[f(x),y],~\forall~x,y \in \mathscr{G}.
    \end{eqnarray*}
    Then $\varphi$ is a super-biderivation.
\end{lemm}

\section{Biderivations on $\mathcal {W}^q$, $\mathcal {V}^q$ and $\mathcal {L}^q$}
 In this section, we first compute biderivations on the $q$-deformed $W(2,2)$ algebra $\mathcal {W}^q$, $q$-deformed Witt algebra $\mathcal {V}^q$ and $q$-deformed Witt superalgebra $\mathcal {L}^q$. Then we characterize their linear commuting maps by using the corresponding results of biderivations.

\subsection{Biderivations on $\mathcal {W}^q$}


First we introduce a special biderivation on $\mathcal {W}^q$. Let $\varphi_0$ be a bilinear map on  $\mathcal {W}^q$ defined by
\begin{align*}
    \varphi_0(L_m,L_n)= & ~[n-m]W_{m+n}, \
    \varphi_0(L_m,W_n)=\varphi_0(W_m,L_n)=\varphi_0(W_m,W_n)=0,~\forall~m,n\in \Z.
\end{align*}
It is easy to check that $\varphi_0$ is a skew-symmetric biderivation on $\mathcal {W}^q$.

\begin{theo} \label{th1}
   Every biderivation on $\mathcal {W}^q$ is of the form $\lambda\varphi_{ad}+\mu \varphi_0$, for some $\lambda, \mu \in \C$.
\end{theo}
\begin{proof}
    Let $\varphi$ be a homogeneous biderivation on $\mathcal {W}^q$ of degree $s$. We can assume that
    \begin{align}\label{assumption}
        \left\{
        \begin{array}{ll}
            \varphi(L_m,L_n)=~a_{s,m,n}^1L_{m+n+s}+a_{s,m,n}^2W_{m+n+s}, \vspace{1.5ex} \\
            \varphi(L_m,W_n)=~b_{s,m,n}^1L_{m+n+s}+b_{s,m,n}^2W_{m+n+s},\vspace{1.5ex}  \\
            \varphi(W_m,L_n)=~c_{s,m,n}^1L_{m+n+s}+c_{s,m,n}^2W_{m+n+s},\vspace{1.5ex}  \\
            \varphi(W_m,W_n)=~d_{s,m,n}^1L_{m+n+s}+d_{s,m,n}^2W_{m+n+s},
        \end{array}
        \right.
    \end{align}
    where $a_{s,m,n}^i,\  b_{s,m,n}^i,\ c_{s,m,n}^i,\ d_{s,m,n}^i\in\C$, $i=1,2$. When there is no ambiguity with the degree $s$, the coefficients $a_{s,m,n}^i,$ $ b_{s,m,n}^i,$ $c_{s,m,n}^i$, $d_{s,m,n}^i$ are simply denoted by $a_{m,n}^i, b_{m,n}^i,c_{m,n}^i,d_{m,n}^i$, respectively.

    Taking the triple $(x,y,z)$ to be $(L_m,L_n,L_p)$,  $(L_m,L_n,W_p)$, $(L_m,W_n,L_p)$ and $(L_m,W_n,W_p)$ in \eqref{vf1}, respectively, gives
    \begin{align}
        \varphi([L_m,L_n],\alpha(L_p))= & ~[\varphi(L_m,L_p),\alpha(L_n)]+[\alpha(L_m),\varphi(L_n,L_p)],\label{vf1(1)} \\
        \varphi([L_m,L_n],\alpha(W_p))= & ~[\varphi(L_m,W_p),\alpha(L_n)]+[\alpha(L_m),\varphi(L_n,W_p)],\label{vf1(2)} \\
        \varphi([L_m,W_n],\alpha(L_p))= & ~[\varphi(L_m,L_p),\alpha(W_n)]+[\alpha(L_m),\varphi(W_n,L_p)],\label{vf1(3)} \\
        \varphi([L_m,W_n],\alpha(W_p))= & ~[\varphi(L_m,W_p),\alpha(W_n)]+[\alpha(L_m),\varphi(W_n,W_p)].\label{vf1(4)}
    \end{align}
    By \eqref{H-L1}, \eqref{alpha1}, \eqref{assumption}-\eqref{vf1(4)}, we obtain the following identities ($i=1,2$):
    \begin{eqnarray}
        &&(q^n+q^{-n})[m+p+s-n]a_{m,p}^i-(q^m+q^{-m})[n+p+s-m]a_{n,p}^i+(q^p+q^{-p})[n-m]a_{m+n,p}^i=0,\label{vf-1}\\
        &&(q^n+q^{-n})[m+p+s-n]b_{m,p}^i-(q^m+q^{-m})[n+p+s-m]b_{n,p}^i+(q^p+q^{-p})[n-m]b_{m+n,p}^i=0,\label{vf-2}\\
        &&(q^m+q^{-m})[n+p+s-m]c_{n,p}^1-(q^p+q^{-p})[n-m]c_{m+n,p}^1=0,\label{vf-34}  \\
        &&(q^n+q^{-n})[m+p+s-n]a_{m,p}^1-(q^m+q^{-m})[n+p+s-m]c_{n,p}^2+(q^p+q^{-p})[n-m]c_{m+n,p}^2=0,\label{vf-34-1}\\
        &&(q^m+q^{-m})[n+p+s-m]d_{n,p}^1-(q^p+q^{-p})[n-m]d_{m+n,p}^1=0,\label{vf-56}\\
        &&(q^n+q^{-n})[m+p+s-n]b_{m,p}^1-(q^m+q^{-m})[n+p+s-m]d_{n,p}^2+(q^p+q^{-p})[n-m]d_{m+n,p}^2=0.\label{vf-56-1}
    \end{eqnarray}

    Similarly, we take the triple $(x,y,z)$ to be $(L_m,L_n,L_p)$,  $(L_m,L_n,W_p)$,  $(W_m,L_n,L_p)$,  $(W_m,L_n,W_p)$ in \eqref{vf2} and obtain
    \begin{align}
        \varphi(\alpha(L_m),[L_n,L_p])= & ~[\varphi(L_m,L_n),\alpha(L_p)]+[\alpha(L_n),\varphi(L_m,L_p)],\label{vf2(1)} \\
        \varphi(\alpha(L_m),[L_n,W_p])= & ~[\varphi(L_m,L_n),\alpha(W_p)]+[\alpha(L_n),\varphi(L_m,W_p)],\label{vf2(2)} \\
        \varphi(\alpha(W_m),[L_n,L_p])= & ~[\varphi(W_m,L_n),\alpha(L_p)]+[\alpha(L_n),\varphi(W_m,L_p)],\label{vf2(3)} \\
        \varphi(\alpha(W_m),[L_n,W_p])= & ~[\varphi(W_m,L_n),\alpha(W_p)]+[\alpha(L_n),\varphi(W_m,W_p)].\label{vf2(4)}
    \end{align}
    By \eqref{H-L1}, \eqref{alpha1} and \eqref{assumption}, \eqref{vf2(1)}-\eqref{vf2(4)} reduce to the following equations ($i=1,2$):
    \begin{eqnarray}
        &&(q^n+q^{-n})[m+p+s-n]a_{m,p}^i-(q^m+q^{-m})[p-n]a_{m,n+p}^i-(q^p+q^{-p})[m+n+s-p]a_{m,n}^i=0,\label{vf-7}\\
        && (q^n+q^{-n})[m+p+s-n]b_{m,p}^1-(q^m+q^{-m})[p-n]b_{m,n+p}^1=0,\label{vf-89}\\
        && (q^n+q^{-n})[m+p+s-n]b_{m,p}^2-(q^m+q^{-m})[p-n]b_{m,n+p}^2-(q^p+q^{-p})[m+n+s-p]a_{m,n}^1=0,\label{vf-89-1}\\
        &&  (q^n+q^{-n})[m+p+s-n]c_{m,p}^i-(q^m+q^{-m})[p-n]c_{m,n+p}^i-(q^p+q^{-p})[m+n+s-p]c_{m,n}^i=0,\label{vf-10}\\
        &&  (q^n+q^{-n})[m+p+s-n]d_{m,p}^1-(q^m+q^{-m})[p-n]d_{m,n+p}^1=0,\label{vf-1112}\\
        &&  (q^n+q^{-n})[m+p+s-n]d_{m,p}^2-(q^m+q^{-m})[p-n]d_{m,n+p}^2-(q^p+q^{-p})[m+n+s-p]c_{m,n}^1=0.\label{vf-1112-1}
    \end{eqnarray}

    In the following we consider two cases:

    \begin{case}{\rm $s\neq 0$

            In \eqref{vf-7} we set $m=p=0$, $n=s$ and obtain $2[s]a^i_{0,s}-2[2s]a^i_{0,s}=0$. Hence $a^i_{0,s}=0$ since $s\neq0$. With this and taking $n=0,$ $p=s$ in \eqref{vf-1}, we have $2[m+2s]a^i_{m,s}+(q^s+q^{-s})[-m]a^i_{m,s}=0$. Thus $a^i_{m,s}=0$ for all $m\in\Z$. Taking this into account and setting $n=s$, $p=0$ in \eqref{vf-7}, we have $(q^s+q^{-s})[m]a^i_{m,0}=0$. This gives
            \begin{eqnarray}\label{1}
                a^i_{m,0}=0,\ \  \forall\, m\in\Z\setminus\{0\}.
            \end{eqnarray}
            Now we set $m=p=0$ in \eqref{vf-1} and obtain
            \begin{eqnarray}\label{2}
                (q^n+q^{-n})[s-n]a^i_{0,0}-2[n+s]a^i_{n,0}+2[n]a^i_{n,0}=0.
            \end{eqnarray}
            We can always find an $n\in\Z$ such that $n\neq 0$ and $n\neq s$. Then it follows from \eqref{1} and \eqref{2} that $a^i_{0,0}=0$. Thus we have proved
                $a^i_{m,0}=0,\ \forall\, m\in\Z.$
            With this we take $p=0$ in \eqref{vf-7} and get
            \begin{eqnarray*}
                (q^m+q^{-m})[n]a^i_{m,n}=2[m+n+s]a^i_{m,n}.
            \end{eqnarray*}
            It follows that
            \begin{eqnarray}\label{4}
                a^i_{m,n}=0,\ \  \forall\, m,n\in\Z; \ \, i=1,2.
            \end{eqnarray}


            With $a_{m,n}^1=0$, and comparing \eqref{vf-89} with \eqref{vf-89-1}, we have
            \begin{eqnarray}
                (q^n+q^{-n})[m+p+s-n]b_{m,p}^i-(q^m+q^{-m})[p-n]b_{m,n+p}^i=0,~i=1,2.\label{89}
            \end{eqnarray}
            Taking $n=p$ in \eqref{89}, we get $(q^p+q^{-p})[m+s]b_{m,p}^i=0$. This gives $b_{m,p}^i=0$ for all $m\neq -s$ and $p\in\Z$.
            Letting $n=0$, $m=-s$ in \eqref{vf-2}, together with $b_{0,p}^i=0$, we have

            \begin{eqnarray*}
                2[p]b_{-s,p}^i+(q^p+q^{-p})[s]b_{-s,p}^i=0.
            \end{eqnarray*}
            It follows $b^i_{-s,p}=0$, and thus
            \begin{eqnarray}\label{b}
                b^i_{m,p}=0,\ \  \forall\, m,p\in\Z; \ i=1,2.
            \end{eqnarray}


            By $a_{m,p}^1=0$ and combining \eqref{vf-34} with \eqref{vf-34-1}, we obtain
            \begin{eqnarray}
                (q^m+q^{-m})[n+p+s-m]c_{n,p}^i-(q^p+q^{-p})[n-m]c_{m+n,p}^i=0, \ \, i=1,2. \label{34}
            \end{eqnarray}
            Letting $n=m$ in \eqref{34} gives $ (q^m+q^{-m})[p+s]c_{m,p}^i=0$. Hence $c_{m,p}^i=0$ for all $p\neq -s$, $m\in\Z$. With $c_{m,0}^i=0$, we set $n=0$, $p=-s$ in \eqref{vf-10} and obtain
                $2[m]c_{m,-s}^i=(q^m+q^{-m})[-s]c_{m,-s}^i.$
            It follows $c_{m,-s}^i=0$. Thus
            \begin{eqnarray}\label{c}
                c_{m,p}^i=0,~\forall~ m,p \in \Z; \ i=1,2.
            \end{eqnarray}

            By the same discussions to $b_{m,p}^i$, one can deduce
            \begin{eqnarray}\label{d}
                d_{m,p}^i=0,~\forall~ m,p \in \Z; \ i=1,2.
            \end{eqnarray}
 By \eqref{4}, \eqref{b}, \eqref{c} and \eqref{d}, we see that the subspaces $\Der_s(\mathcal {W}^q,\mathcal {W}^q)$ vanish when $s\neq 0$.
        }\end{case}

    \begin{case}{\rm $s=0$

            Letting $n=p$ in \eqref{vf-1}, we get
            \begin{eqnarray}
                (q^n+q^{-n})[m]a_{m,n}^i-(q^m+q^{-m})[2n-m]a_{n,n}^i+(q^n+q^{-n})[n-m]a_{m+n,n}^i=0.\label{1:n=p}
            \end{eqnarray}
             Setting $m=0$ in \eqref{1:n=p} gives $\big( (q^n+q^{-n})[n]-2[2n]\big)a_{n,n}^i=0$.  This implies $a_{n,n}^i=0$ for all $n\in\Z\setminus \{0\}$. Taking $n=1$ in \eqref{1:n=p} and with $a_{1,1}^i=0$, we have
            \begin{eqnarray*}
                (q+q^{-1})[m]a_{m,1}^i+(q+q^{-1})[1-m]a_{m+1,1}^i=0,
            \end{eqnarray*}
            which implies
            \begin{eqnarray}\label{a-1}
                \left\{
                \begin{array}{ll}
                    a_{m,1}^i=[1-m]a_{0,1}^i, ~m\leq 1,\vspace{2ex} \\
                    a_{m,1}^i=[m-1]a_{2,1}^i, ~m\geq 2.
                \end{array}
                \right.
            \end{eqnarray}
            Taking $m=2$, $n=0$ and $p=1$ in \eqref{vf-1}, we get
            \begin{eqnarray*}
                2[3]a_{2,1}^i+(q^2+q^{-2})a_{0,1}^i-(q+q^{-1})[2]a_{2,1}^i=0,
            \end{eqnarray*}
            which gives $a_{2,1}^i=-a_{0,1}^i$. Combining with \eqref{a-1}, we obtain
            \begin{eqnarray}\label{a-2}
            a_{m,1}^i=[1-m]a_{0,1}^i, ~\forall~ m\in \Z.
            \end{eqnarray}

          Taking $n=0$, $p=1$ in \eqref{vf-7}, we have
            \begin{eqnarray*}
                2[m+1]a_{m,1}^i-(q^m+q^{-m})a_{m,1}^i-(q+q^{-1})[m-1]a_{m,0}^i=0.
            \end{eqnarray*}
          This, together with \eqref{a-2}, gives
            \begin{eqnarray}\label{a-4}
                a_{m,0}^i=\frac{2[m+1]-(q^m+q^{-m})}{(q+q^{-1})[m-1]}a_{m,1}^i=[-m]a_{0,1}^i,~~m\neq 1.
            \end{eqnarray}
             Setting $m=2$, $n=-1$, $p=0$ in \eqref{vf-1} gives
            \begin{eqnarray}\label{a-3}
                (q+q^{-1})[3]a_{2,0}^i+(q^2+q^{-2})[3]a_{-1,0}^i-2[3]a_{1,0}^i=0.
            \end{eqnarray}
            Inserting $a_{2,0}^i=[-2]a_{0,1}^i$ and $a_{-1,0}^i=a_{0,1}^i$ \big(see \eqref{a-4}\big) into \eqref{a-3}, we get $a_{1,0}^i=[-1]a_{0,1}^i$. Combining with \eqref{a-4}, we obtain
           \begin{eqnarray}\label{a-5}
            a_{m,0}^i=[-m]a_{0,1}^i,~\forall ~m\in \Z.
             \end{eqnarray}
            Taking $p=0$ in \eqref{vf-7} gives
            \begin{eqnarray*}
                (q^n+q^{-n})[m-n]a_{m,0}^i+(q^m+q^{-m})[n]a_{m,n}^i-2[m+n]a_{m,n}^i=0,
            \end{eqnarray*}
           from which, together with \eqref{a-5}, we can deduce that
            \begin{eqnarray}\label{a-6}
                a_{m,n}^i=\frac{(q^n+q^{-n})[m-n]}{2[m+n]-(q^m+q^{-m})[n]}a_{m,0}^i=\frac{[m-n]}{[m]}a_{m,0}^i=[n-m]a_{0,1}^i,~m\neq 0.
            \end{eqnarray}
            Letting $m=1$, $n=-1$ in \eqref{vf-1} gives
            \begin{eqnarray}\label{a-7}
                (q+q^{-1})[2+p]a_{1,p}^i-(q+q^{-1})[p-2]a_{-1,p}^i+(q^p+q^{-p})[-2]a_{0,p}^i=0.
            \end{eqnarray}
            By \eqref{a-6}, we have $a_{1,p}^i=[p-1]a_{0,1}^i$ and $~a_{-1,p}^i=[p+1]a_{0,1}^i$. Substituting this into \eqref{a-7} gives
            $a_{0,p}^i=[p]a_{0,1}^i,$ for all $p\in \Z$. Combining this with \eqref{a-6}, we obtain
            \begin{eqnarray}\label{aa}
                a_{m,n}^i=[n-m]a_{0,1}^i,~\forall~ m,n\in \Z;~i=1,2.
            \end{eqnarray}

            Setting $n=0$ in \eqref{vf-89}, we obtain
            \begin{eqnarray*}
                2[m+p]b_{m,p}^1-(q^m+q^{-m})[p]b_{m,p}^1=(q^p+q^{-p})[m]b_{m,p}^1=0,
            \end{eqnarray*}
            which gives $b_{m,p}^1=0$ for all $m\neq 0$, $p\in\Z$. In particular, $b_{1,p}^1=b_{-1,p}^1=0$. With this and
            taking $m=1$, $n=-1$ in \eqref{vf-2}, we have $b_{0,p}^1=0$ for all $p\in\Z$. Hence
             \begin{eqnarray}\label{bb}
             b_{m,p}^1=0,~\forall ~m,p\in \Z.
              \end{eqnarray}

            In \eqref{vf-89-1} we set $p=n$ and get
            \begin{eqnarray*}
                (q^n+q^{-n})[m]b_{m,n}^2-(q^n+q^{-n})[m]a_{m,n}^1=0.
            \end{eqnarray*}
           This, together with \eqref{aa}, gives
             \begin{eqnarray}\label{b-2}
            b_{m,n}^2=a_{m,n}^1=[n-m]a_{0,1}^1, \ \forall \, n\in\Z, \,m\in\Z\setminus\{0\}.
             \end{eqnarray}
            Taking $m=1$, $n=-1$ in \eqref{vf-2} gives
            \begin{eqnarray}\label{b-1}
                (q+q^{-1})[p+2]b_{1,p}^2-(q+q^{-1})[p-2]b_{-1,p}^2+(q^p+q^{-p})[-2]b_{0,p}^2=0.
            \end{eqnarray}
            Inserting $b_{1,p}^2=[p-1]a_{0,1}^1$ and  $b_{-1,p}^2=[p+1]a_{0,1}^1$ into \eqref{b-1}, we have
            \begin{eqnarray*}
                b_{0,p}^2=\frac{[p+2]}{q^p+q^{-p}}b_{1,p}^2-\frac{[p-2]}{q^p+q^{-p}}b_{-1,p}^2=[p]a_{0,1}^1.
            \end{eqnarray*}
          Combining this with \eqref{b-2}, we obtain
           \begin{eqnarray}\label{b-3}
           b_{m,p}^2=[p-m]a_{0,1}^1,~\forall~m,p\in \Z.
           \end{eqnarray}

            Letting $m=0$ in \eqref{vf-34}, we obtain $ 2[n+p]c_{n,p}^1-(q^p+q^{-p})[n]c_{n,p}^1=(q^n+q^{-n})[p]c_{n,p}^1=0$. Hence $c_{n,p}^1=0,$ $\forall ~p\neq 0$ and $n\in\Z$. Taking this into account and setting $p=1$, $n=-1$ in \eqref{vf-10},  we get  $c_{m,0}^1=0$ for all $m\in\Z$. Thus
            \begin{eqnarray}\label{c-1}
            c_{m,p}^1=0,~\forall~ m,p\in \Z.
            \end{eqnarray}

            Setting $m=n$ in \eqref{vf-34-1}, we obtain
            \begin{eqnarray*}
                (q^n+q^{-n})[p]a_{n,p}^1-(q^n+q^{-n})[p]c_{n,p}^2=0.
            \end{eqnarray*}
             This, together with \eqref{aa}, yields
            \begin{eqnarray}\label{c-2}
            c_{n,p}^2=a_{n,p}^1=[p-n]a_{0,1}^1, \ \forall \, n\in\Z, \,p\in\Z\setminus\{0\}.
             \end{eqnarray}
            Taking $p=1$, $n=-1$ in \eqref{vf-10} gives
            \begin{eqnarray*}
                (q+q^{-1})[m+2]c_{m,1}^2-(q^m+q^{-m})[2]c_{m,0}^2-(q+q^{-1})[m-2]c_{m,-1}^2=0.
            \end{eqnarray*}
            Plugging $c_{m,1}^2=[1-m]a_{0,1}^1$ and $c_{m,-1}^2=[-1-m]a_{0,1}^1$ into the above equation, we have
            \begin{eqnarray*}
                c_{m,0}^2=\frac{[m+2]}{q^m+q^{-m}}c_{m,1}^2-\frac{[m-2]}{q^m+q^{-m}}c_{m,-1}^2=[-m]a_{0,1}^1.
            \end{eqnarray*}
            Combining this with \eqref{c-2}, we get
            \begin{eqnarray}\label{c-3}
            c_{m,p}^2=[p-m]a_{0,1}^1,~\forall~ m,p\in \Z.
            \end{eqnarray}

            Since $c_{m,n}^1=0$, \eqref{vf-1112} and \eqref{vf-1112-1} reduce to
             \begin{eqnarray*}
                (q^n+q^{-n})[m+p-n]d_{m,p}^i-(q^m+q^{-m})[p-n]d_{m,n+p}^i=0,\ \,i=1,2.\label{1112}
            \end{eqnarray*}
            Comparing this with \eqref{vf-89} and using the same discussions to $b^1_{m,p}$, we obtain
            \begin{eqnarray}\label{dd}
           d_{m,p}^i=0,~\forall~m,p\in \Z;\ \,i=1,2.
           \end{eqnarray}

           Set $a^1_{0,1}=\lambda$ and $a^2_{0,1}=\mu$. By \eqref{assumption}, \eqref{aa}, \eqref{bb}, \eqref{b-3}, \eqref{c-1}, \eqref{c-3} and \eqref{dd}, we conclude that $\varphi=\lambda\varphi_{ad}+\mu \varphi_0$ in the case of $s=0$.
        }\end{case}

    Combining the two cases above, the proof is complete.
\end{proof}

\subsection{Biderivations on $\mathcal {V}^q$}


Note that the $q$-deformed  Witt algebra $\mathcal {V}^q$ has a $\Z$-grading, i.e.,
$\mathcal {V}^q=\bigoplus_{n\in\mathbb{Z}}\mathcal {V}^q_n, \ \ \ \mathcal {V}^q_n={\rm span}_{\mathbb{C}}\{L_n\}.$

\begin{theo} \label{th2}
    Every biderivation of $\mathcal {V}^q$ is inner.
\end{theo}

\begin{proof}
    Let $\varphi$ be a biderivation on $\mathcal{V}^q$ of degree $s$. Assume that
    \begin{eqnarray}\label{vv}
        \varphi(L_m,L_n)=a_{m,n}L_{m+n+s},~~where ~a_{m,n}\in\C.
    \end{eqnarray}
    Taking the triple $(x,y,z)$ to be $(L_m,L_n,L_p)$ in \eqref{vf1} and \eqref{vf2} respectively gives
    \begin{align*}
        \varphi([L_m,L_n],\alpha(L_p))= & ~[\varphi(L_m,L_p),\alpha(L_n)]+[\alpha(L_m),\varphi(L_n,L_p)], \\
        \varphi(\alpha(L_m),[L_n,L_p])= & ~[\varphi(L_m,L_n),\alpha(L_p)]+[\alpha(L_n),\varphi(L_m,L_p)].
    \end{align*}
    Inserting \eqref{vv} into the above equations respectively, we have
    \begin{eqnarray}
        (1+q^n)(\{n\}-\{m+p+s\})a_{m,p}+(1+q^m)(\{n+p+s\}-\{m\})a_{n,p}-(1+q^p)(\{n\}-\{m\})a_{m+n,p}=0,\label{witt-1}\\
        (1+q^n)(\{m+p+s\}-\{n\})a_{m,p}-(1+q^m)(\{p\}-\{n\})a_{m,n+p}+(1+q^p)(\{p\}-\{m+n+s\})a_{m,n}=0.\label{witt-2}
    \end{eqnarray}

    We first consider the case $s\neq 0$. In \eqref{witt-2} we set $m=p=0$, $n=s$ and obtain $2\{s\}a_{0,s}-2\{2s\}a_{0,s}=0$. It implies $a_{0,s}=0$ since $s\neq0$. With this and taking $n=0$, $p=s$ in \eqref{witt-1}, we get $-2\{m+2s\}a_{m,s}+(1+q^{s})\{m\}a_{m,s}=0$. Hence $a_{m,s}=0$ for all $m\in\Z$. Taking this into account and setting $n=s$, $p=0$ in \eqref{witt-2}, we have $(1+q^{s})(\{m+s\}-\{s\})a_{m,0}=0$. This gives

    \begin{eqnarray}\label{w1}
        a_{m,0}=0,\ \  \forall\, m\in\Z\setminus\{0\}.
    \end{eqnarray}
    Then we set $m=p=0$ in \eqref{witt-1} and obtain
    \begin{eqnarray}\label{w2}
        (1+q^{n})(\{n\}-\{s\})a_{0,0}+2\{n+s\}a_{n,0}-2\{n\}a_{n,0}=0.
    \end{eqnarray}
    We can always find an $n\in\Z$ such that $n\neq 0$ and $n\neq s$. Then it follows from \eqref{w1} and \eqref{w2} that $a_{0,0}=0$. Thus we have proved that
    \begin{eqnarray}
        a_{m,0}=0,\ \  \forall\, m\in\Z.
    \end{eqnarray}
    With this, we take $p=0$ in \eqref{witt-2} and get
    \begin{eqnarray*}
        (1+q^{m})\{n\}a_{m,n}=2\{m+n+s\}a_{m,n}.
    \end{eqnarray*}
    It follows that
    \begin{eqnarray}
        a_{m,n}=0,\ \  \forall\, m,n\in\Z.
    \end{eqnarray}

    Now we discuss the case $s=0$. Taking $m=p=0$ and $n=1$ in \eqref{witt-1} gives $a_{0,0}=0$. Letting $p=n$ in \eqref{witt-1}, we get
    \begin{eqnarray}
        (1+q^n)(\{n\}-\{m+n\})a_{m,n}+(1+q^m)(\{2n\}-\{m\})a_{n,n}-(1+q^n)(\{n\}-\{m\})a_{m+n,n}=0.\label{witt-1:p=n}
    \end{eqnarray}
     We set $m=0$ in \eqref{witt-1:p=n} and obtain $ \big(2\{2n\}-(1+q^n)\{n\}\big)a_{n,n}=\{2n\}a_{n,n}=0$. Hence $a_{n,n}=0$ for $n\in\Z$. Then taking $n=1$ in \eqref{witt-1:p=n} and with $a_{1,1}=0$, we have
    \begin{eqnarray*}
        (1+q)(\{1\}-\{m+1\})a_{m,1}-(1+q)(\{1\}-\{m\})a_{m+1,1}=0,
    \end{eqnarray*}
    which implies
    \begin{eqnarray}\label{aaa}
        \left\{
        \begin{array}{ll}
            a_{m,1}=(\{1\}-\{m\})a_{0,1}, ~m\leq 1,\vspace{2ex} \\
            a_{m,1}=\frac{\{1\}-\{m\}}{\{1\}-\{2\}}a_{2,1}, ~m\geq 2.
        \end{array}
        \right.
    \end{eqnarray}
    Letting $m=2$, $n=0$, $p=1$ in \eqref{witt-1}, we get
    \begin{eqnarray*}
        -2\{3\}a_{2,1}+(1+q^2)(\{1\}-\{2\})a_{0,1}+(1+q)\{2\}a_{2,1}=0.
    \end{eqnarray*}
    It follows that $a_{2,1}=\frac{(1+q^2)(\{1\}-\{2\})}{2\{3\}-(1+q)\{2\}}a_{0,1}=(\{1\}-\{2\})a_{0,1}$. Hence, by \eqref{aaa}, we obtain
    \begin{eqnarray}\label{aaa-1}
    a_{m,1}=(\{1\}-\{m\})a_{0,1}, ~\forall~ m\in \Z.
    \end{eqnarray}

   Setting $n=0$ and $p=1$ in \eqref{witt-2} gives
    \begin{eqnarray*}
        2\{m+1\}a_{m,1}-(1+q^m)a_{m,1}+(1+q)(1-\{m\})a_{m,0}=0,
    \end{eqnarray*}
    and hence
    \begin{eqnarray}\label{aaa-2}
        a_{m,0}=\frac{(1+q^m)-2\{m+1\}}{(1+q)(1-\{m\})}a_{m,1}=-\{m\}a_{0,1},~m\neq 1.
    \end{eqnarray}
     Setting $m=1$, $n=-1$, $p=0$ in \eqref{witt-1} and by using $a_{0,0}=0$, we have
    \begin{eqnarray*}
        (1+q^{-1})(\{-1\}-1)a_{1,0}+(1+q)(\{-1\}-\{1\})a_{-1,0}=0.
    \end{eqnarray*}
   This, together with $a_{-1,0}=-\{-1\}a_{0,1}=q^{-1}a_{0,1}$, yields
   $a_{1,0}=-\frac{1+q}{1+q^{-1}}a_{-1,0}=-a_{0,1}$. Hence, by \eqref{aaa-2}, we obtain
   \begin{eqnarray}\label{aaa-3}
   a_{m,0}=-\{m\}a_{0,1},~\forall ~m\in \Z.
   \end{eqnarray}
 Taking $p=0$ in \eqref{witt-2} gives
    \begin{eqnarray}\label{aaa-4}
        (1+q^n)(\{m\}-\{n\})a_{m,0}+(1+q^m)\{n\}a_{m,n}-2\{m+n\}a_{m,n}=0.
    \end{eqnarray}
    Substituting \eqref{aaa-3} into \eqref{aaa-4}, we can deduce that
    \begin{eqnarray}\label{aaa-5}
        a_{m,n}=\frac{(1+q^n)(\{m\}-\{n\})}{2\{m+n\}-(1+q^m)\{n\}}a_{m,0}=(\{n\}-\{m\})a_{0,1},~\forall~m\neq 0, \ n\in\Z.
    \end{eqnarray}
    Letting $m=1$, $n=-1$ in \eqref{witt-1} gives
    \begin{eqnarray}\label{aaa-6}
        (1+q^{-1})(\{-1\}-\{p+1\})a_{1,p}+(1+q)(\{p-1\}-\{1\})a_{-1,p}-(1+q^p)(\{-1\}-\{1\})a_{0,p}=0.
    \end{eqnarray}
    By \eqref{aaa-5}, we have $a_{1,p}=(\{p\}-\{1\})a_{0,1}$ and $~a_{-1,p}=(\{p\}-\{-1\})a_{0,1}$. Inserting this into \eqref{aaa-6}, we obtain
    \[
        a_{0,p}=\frac{(1+q^{-1})(\{-1\}-\{p+1\})}{(1+q^p)(\{-1\}-\{1\})}a_{1,p}+\frac{(1+q)(\{p-1\}-\{1\})}{(1+q^p)(\{-1\}-\{1\})}a_{-1,p}=\{p\}a_{0,1},~\forall~ p\in \Z.
    \]
    Combining this with \eqref{aaa-5}, we get
    \begin{eqnarray*}
        a_{m,n}=(\{n\}-\{m\})a_{0,1},~\forall~ m,n\in \Z.
    \end{eqnarray*}
    Setting $a_{0,1}=\lambda$ and by the discussions above, we have
    \begin{eqnarray*}
        \varphi(L_m,L_n)=\lambda (\{n\}-\{m\}) L_{m+n}=\lambda[L_m,L_n],~\forall~m,n\in\Z.
    \end{eqnarray*}
    This proves the theorem.
\end{proof}

\subsection{Super-biderivations on $\mathcal{L}^q$}

It is easy to see that $q$-deformed Witt superalgebra $\mathcal{L}^q$ has a $\Z$-grading, i.e.,
$\mathcal {L}^q=\bigoplus_{n\in\mathbb{Z}}\mathcal {L}^q_n$ with $ \mathcal {L}^q_n={\rm span}_{\mathbb{C}}\{L_n,G_n\}.$
We have
$${\rm BDer}(\mathcal{L}^q,\mathcal{L}^q) = {\rm BDer}_{\bar 0}(\mathcal{L}^q,\mathcal{L}^q)\oplus {\rm BDer}_{\bar 1}(\mathcal{L}^q,\mathcal{L}^q).$$
Denote by ${\rm BDer}_{\bar 0,s}(\mathcal{L}^q,\mathcal{L}^q)$ (resp. ${\rm BDer}_{\bar1,s}(\mathcal{L}^q,\mathcal{L}^q)$) the subspace of ${\rm BDer}_{\bar 0}(\mathcal{L}^q,\mathcal{L}^q)$ (resp. ${\rm BDer}_{\bar 1}(\mathcal{L}^q,\mathcal{L}^q)$) of degree $s$. That is given by even (resp. odd) super-biderivations $\varphi$
of degree $s$, i.e., for all homogeneous elements $x_m, x_n\in \mathcal {L}^q$  of degree $m$ and $n$ respectively,  $\varphi(x_m, x_n)$ is of degree $m + n + s$.

Define a bilinear map $\varphi_{-1}$ of $\mathcal{L}^q$ on the generators by
\begin{eqnarray}\label{verfi-1}
\varphi_{-1}(L_m,L_n)=(\{n\}-\{m\})G_{m+n-1}, \ \varphi_{-1}(L_m,G_n)=\varphi_{-1}(G_m,L_n)=\varphi_{-1}(G_m,G_n)=0, \ \,\forall\ \,m,n\in\Z.
\end{eqnarray}
It is easy to check that $\varphi_{-1}$ is an odd super-biderivation of degree $-1$.

Our main theorem in this subsection is the following:

\begin{theo} \label{th3}
    \begin{itemize}
        \item[(1)] Any even super-biderivation of the $q$-deformed Witt superalgebra is inner.
        \item[(2)] Any odd super-biderivation of the $q$-deformed Witt superalgebra is of the form $\lambda\varphi_{-1}$, for some $\lambda\in\mathbb{C}$.
    \end{itemize}
\end{theo}

In the sequel, we proceed to prove this theorem by computing the even and odd
subspaces of super-biderivations of $\mathcal{L}^q$, respectively.

\subsubsection{Even super-biderivations}

 Let $\varphi$ be an even super-biderivation on $\mathcal {L}^q$ of degree $s$. We can assume that
\begin{align*}
    \varphi(L_m,L_n)=a_{m,n}L_{m+n+s},~
    \varphi(L_m,G_n)=b_{m,n}G_{m+n+s},~
    \varphi(G_m,L_n)=c_{m,n}G_{m+n+s},~
    \varphi(G_m,G_n)=d_{m,n}L_{m+n+s},
\end{align*}
where $a_{m,n}, b_{m,n}, c_{m,n}, d_{m,n}\in\C$.

Taking the triple $(x,y,z)$ to be $(L_m,L_n,G_p)$, $(G_m,L_n,L_p)$ and $(G_m,L_n,G_p)$ in \eqref{svf1} and \eqref{svf2}, respectively, we obtain the following identities:


\begin{eqnarray}
    \nonumber   && (1+q^n)(\{n\}-\{m+p+s+1\})b_{m,p} +(1+q^m)(\{n+p+s+1\}-\{m\})b_{n,p}                      \\
                                                     &&~~~~~~~~~~~~~~~~~~~~~~~~~~~~~~~~~~~~~~~~~~~~~~~~~~~~~~~~~~~ -(1+q^{p+1})(\{n\}-\{m\})b_{m+n,p}=0,\label{even-1}     \\
    \nonumber   && (1+q^n)(\{m+p+s+1\}-\{n\})b_{m,p}-(1+q^m)(\{p+1\}-\{n\})b_{m,n+p}                        \\
                                                     &&~~~~~~~~~~~~~~~~~~~~~~~~~~~~~~~~~~~~~~~~~~~~~~~~~~~~~~~~~~~+(1+q^{p+1})(\{p+1\}-\{m+n+s\})a_{m,n}=0,\label{even-2}\\ \nonumber     &&(1+q^n)(\{n\}-\{m+p+s+1\})c_{m,p}+(1+q^{m+1})(\{n+p+s\}-\{m+1\})a_{n,p}                  \\
                                                     &&~~~~~~~~~~~~~~~~~~~~~~~~~~~~~~~~~~~~~~~~~~~~~~~~~~~~~~~~~~~ -(1+q^p)(\{n\}-\{m+1\})c_{m+n,p}=0,\label{even-3}       \\
    \nonumber     && (1+q^n)(\{m+p+s+1\}-\{n\})c_{m,p}-(1+q^{m+1})(\{p\}-\{n\})c_{m,n+p}                      \\
                                                     &&~~~~~~~~~~~~~~~~~~~~~~~~~~~~~~~~~~~~~~~~~~~~~~~~~~~~~~~~~~~+(1+q^p)(\{p\}-\{m+n+s+1\})c_{m,n}=0,\label{even-4}     \\
    &&(1+q^n)(\{n\}-\{m+p+s\})d_{m,p}-(1+q^{p+1})(\{n\}-\{m+1\})d_{m+n,p}=0,\label{even-5}   \\
    &&(1+q^n)(\{m+p+s\}-\{n\})d_{m,p}-(1+q^{m+1})(\{p+1\}-\{n\})d_{m,n+p}=0.\label{even-6}
\end{eqnarray}




\begin{prop} \label{prop1} If $s\neq 0$, the subspaces ${\rm BDer}_{\bar0,s}(\mathscr{L}^q,\mathscr{L}^q)$ vanish.

\end{prop}
\begin{proof}
     By the proof of Theorem \ref{th2}, we have $a_{m,n}=0$ for all $m,n\in \Z$.
Thus Eq. \eqref{even-2} becomes
    \begin{eqnarray}
        (1+q^n)(\{m+p+s+1\}-\{n\})b_{m,p}-(1+q^m)(\{p+1\}-\{n\})b_{m,n+p}=0.\label{even-2-1}
    \end{eqnarray}
    In \eqref{even-2-1} we set $n=p+1$ and obtain
    \begin{eqnarray*}
        (1+q^{p+1})(\{m+p+s+1\}-\{p+1\}) b_{m,p}=0.
    \end{eqnarray*}
    This implies $b_{m,p}=0$ for $m\neq -s$. Then we set $n=0$, $m=-s$ in \eqref{even-1} and with $b_{0,p}=0$, we have $2\{p+1\}b_{-s,p}=(1+q^{p+1})\{-s\}b_{-s,p}$. This gives $b_{-s,p}=0$. Hence
    $b_{m,p}=0,~\forall~ m,p\in \Z.$

    Since $a_{m,n}=0$,  Eq. \eqref{even-3} reduces to
    \begin{eqnarray}
        (1+q^n)(\{n\}-\{m+p+s+1\})c_{m,p}-(1+q^p)(\{n\}-\{m+1\})c_{m+n,p}=0.\label{even-3-1}
    \end{eqnarray}
    Letting $n=m+1$ in \eqref{even-3-1} gives
    \begin{eqnarray*}
        (1+q^{m+1})(\{m+1\}-\{m+p+s+1\})c_{m,p}=0.
    \end{eqnarray*}
    This implies $c_{m,p}=0$ for $p\neq -s$. With $c_{m,0}=0$, we
    set $p=-s$ and $n=0$ in \eqref{even-4} and obtain $2\{m+1\}c_{m,-s}-(1+q^{m+1})\{-s\})c_{m,-s}=0$, leading to $c_{m,-s}=0$. It follows
    $c_{m,p}=0,~\forall~ m,p\in \Z.$

    Taking $n=m+1$ in \eqref{even-5}, we have
    \begin{eqnarray*}\label{even-5-1}
        (1+q^{m+1})(\{m+1\}-\{m+p+s\})d_{m,p}=0.
    \end{eqnarray*}
    From this equation we can deduce $d_{m,p}=0$ for $p\neq 1-s$. Taking $n=1$, $p=-s$ in \eqref{even-6} and with $d_{m,-s}=0$, we obtain $(1+q^{m+1})(\{1-s\}-\{1\})d_{m,1-s}=0$. This gives $d_{m,1-s}=0$. Hence
    $d_{m,p}=0,~\forall~ m,p\in \Z.$

    The proof is complete.
\end{proof}

Now we discuss the case $s=0$.

\begin{prop} \label{prop2}
    Each even super-biderivation of $\mathscr{L}^q$ of degree 0 is inner.
\end{prop}
\begin{proof}
    By the proof of Theorem \ref{th2}, there exists $\lambda\in\C$ such that $a_{m,n}=(\{n\}-\{m\})\lambda,~\forall~ m,n\in \Z.$
    In \eqref{even-2}, we set $n=0$ and obtain
    \begin{eqnarray*}
        \big( 2\{m+p+1\}-(1+q^m)\{p+1\} \big)b_{m,p}+(1+q^{p+1})(\{p+1\}-\{m\})a_{m,0}=0.
    \end{eqnarray*}
    For $m\neq 0$, we deduce that
    \begin{align}\label{even-5-2}
        b_{m,p}=~\frac{(1+q^{p+1})(\{m\}-\{p+1\})}{ 2\{m+p+1\}-(1+q^m)\{p+1\}}a_{m,n}
        =~(\{p+1\}-\{m\})\lambda.
    \end{align}
    Taking $m=1$, $n=-1$ in \eqref{even-1}, we have
    \begin{eqnarray}\label{even-5-33}
        (1+q^{-1})(\{-1\}-\{p+2\})b_{1,p}+(1+q)(\{p\}-\{1\})b_{-1,p}-(1+q^{p+1})(\{-1\}-\{1\})b_{0,p}=0.
    \end{eqnarray}
    It follows from \eqref{even-5-2} that $b_{1,p}=(\{p+1\}-\{1\})\lambda$ and  $b_{-1,p}=(\{p+1\}-\{-1\})\lambda$.
    Inserting this into \eqref{even-5-33},
     we get $b_{0,p}=\{p+1\}\lambda$. This, together with \eqref{even-5-2}, gives
    \[
        b_{m,p}=(\{p+1\}-\{m\})\lambda,~\forall~ m,n\in \Z.
    \]

    Setting $n=0$ in \eqref{even-3}, we obtain
    \begin{eqnarray*}
        \big( (1+q^p)\{m+1\}-2\{m+p+1\} \big)c_{m,p}+(1+q^{m+1})(\{p\}-\{m+1\})a_{0,p}=0.
    \end{eqnarray*}
    For $p\neq 0$, we have
    \begin{align}\label{even-5-3}
        c_{m,p}=~\frac{(1+q^{m+1})(\{p\}-\{m+1\})}{ 2\{m+p+1\}-(1+q^p)\{m+1\}}a_{0,p}
        =~(\{p\}-\{m+1\})\lambda.
    \end{align}
    Taking  $n=-1$, $p=1$ in \eqref{even-4}, we have
    \begin{eqnarray}
        (1+q^{-1})(\{m+2\}-\{-1\})c_{m,1}-(1+q^{m+1})(\{1\}-\{-1\})c_{m,0}+(1+q^p)(\{1\}-\{m\})c_{m,-1}=0.
    \end{eqnarray}
    Substituting $c_{m,1}=(\{1\}-\{m+1\})\lambda$ and  $c_{m,-1}=(\{-1\}-\{m+1\})\lambda$ into this equation, we get $c_{m,0}=-\{m+1\}\lambda$. This, combining with \eqref{even-5-3}, gives
    \[
        c_{m,p}=(\{p\}-\{m+1\})\lambda,~\forall~ m,p\in \Z.
    \]

    Taking $m=n=0$ in \eqref{even-5}, we obtain $d_{0,p}=0,~\forall~p\in \Z$. Setting $m=0$ in \eqref{even-5} gives
    $d_{n,p}=0,~\forall~ p\in \Z$ and $n\neq 1$. Finally, we set $m=1$ and $n=2$ and obtain $d_{1,p}=0$ for all $p\in \Z$.
   Hence $d_{m,p}=0,~\forall~ m,p\in \Z.$

     We conclude that
    \begin{align*}
        \varphi(L_m,L_n)= & ~\lambda(\{n\}-\{m\})L_{m+n}=\lambda[L_m,L_n],   \\
        \varphi(L_m,G_n)= & ~\lambda(\{n+1\}-\{m\})G_{m+n}=\lambda[L_m,G_n], \\
        \varphi(G_m,L_n)= & ~\lambda(\{n\}-\{m+1\})G_{m+n}=\lambda[G_m,L_n], \\
        \varphi(G_m,G_n)= & ~0.
    \end{align*}
    Hence the proposition is proved.
\end{proof}

Taking Propositions \ref{prop1} and  \ref{prop2} together, we find the announced result in Theorem \ref{th3} (1).

\subsubsection{Odd super-biderivations}

 The results in this subsection shall be the following:
\begin{prop} \label{prop6} The following hold:
    \begin{itemize}
        \item [(1)]  ${\rm BDer}_{\bar 1,s}(\mathcal{L}^q,\mathcal{L}^q)=\{0\}$,  for $s\neq -1$,
        \item [(2)]  ${\rm BDer}_{\bar 1,-1}(\mathcal{L}^q,\mathcal{L}^q)=\langle\varphi_{-1}\rangle$, where $\varphi_{-1}$ is defined by \eqref{verfi-1}.
    \end{itemize}
\end{prop}
\begin{proof}
    Let $\varphi$ be an odd super-biderivation on $\mathcal{L}^q$ of degree $s$. Assume that
    \begin{align*}
        \varphi(L_m,L_n)=  a_{m,n}G_{m+n+s},\
        \varphi(L_m,G_n)=b_{m,n}L_{m+n+s},
        \varphi(G_m,L_n)=  c_{m,n}L_{m+n+s},\
        \varphi(G_m,G_n)=d_{m,n}G_{m+n+s},
    \end{align*}
    where $a_{m,n}, b_{m,n}, c_{m,n}, d_{m,n}\in\C$.

    Taking the triple $(x,y,z)$ to be $(L_m,L_n,L_p)$, $(L_m,L_n,G_p)$, $(G_m,L_n,L_p)$ and $(G_m,L_n,G_p)$ in \eqref{svf1} and \eqref{svf2}, respectively, we have
     \begin{align}
        \nonumber    (1+q^n)(\{n\}-\{m+p+s+1\})a_{m,p}+&~ (1+q^m)(\{n+p+s+1\}-\{m\})a_{n,p}                     \\
                                                         ~~~~~~~~~~~~~~~~~~~~~~~~~~~~~~~~~~~~~~~~~~~~~~~~~~~~ -&~(1+q^{p})(\{n\}-\{m\})a_{m+n,p}=0,\label{odd-1}       \\
        \nonumber     (1+q^n)(\{m+p+s+1\}-\{n\})a_{m,p} -&~(1+q^m)(\{p\}-\{n\})a_{m,n+p}                         \\
                                                         ~~~~~~~~~~~~~~~~~~~~~~~~~~~~~~~~~~~~~~~~~~~~~~~~~~~~ +&~(1+q^{p})(\{p\}-\{m+n+s+1\})a_{m,n}=0,\label{odd-2}   \\
        \nonumber (1+q^n)(\{n\}-\{m+p+s\})b_{m,p}+&~(1+q^m)(\{n+p+s\}-\{m\})b_{n,p}                       \\ ~~~~~~~~~~~~~~~~~~~~~~~~~~~~~~~~~~~~~~~~~~~~~~~~~~~~-&~(1+q^{p+1})(\{n\}-\{m\})b_{m+n,p}=0,\label{odd-3}\\
        (1+q^n)(\{m+p+s\}-\{n\})b_{m,p}-&~(1+q^m)(\{p+1\}-\{n\})b_{m,n+p}=0,\label{odd-4}       \\
         (1+q^n)(\{n\}-\{m+p+s\})c_{m,p}-&~(1+q^p)(\{n\}-\{m+1\})c_{m+n,p}=0,\label{odd-5}       \\
        \nonumber (1+q^n)(\{m+p+s\}-\{n\})c_{m,p}-&~(1+q^{m+1})(\{p\}-\{n\})c_{m,n+p}                     \\ ~~~~~~~~~~~~~~~~~~~~~~~~~~~~~~~~~~~~~~~~~~~~~~~~~~~~+&~(1+q^p)(\{p\}-\{m+n+s\})c_{m,n}=0,\label{odd-6}\\
        \nonumber     (1+q^n)(\{n\}-\{m+p+s+1\})d_{m,p}-&~(1+q^{m+1})(\{n+p+s\}-\{m+1\})b_{n,p}                 \\
                                                         ~~~~~~~~~~~~~~~~~~~~~~~~~~~~~~~~~~~~~~~~~~~~~~~~~~~~ -&~(1+q^{p+1})(\{n\}-\{m+1\})d_{m+n,p}=0,\label{odd-7}   \\
        \nonumber   (1+q^n)(\{m+p+s+1\}-\{n\})d_{m,p}-&~(1+q^{m+1})(\{p+1\}-\{n\})d_{m,n+p}                   \\
                                                         ~~~~~~~~~~~~~~~~~~~~~~~~~~~~~~~~~~~~~~~~~~~~~~~~~~~~ +&~(1+q^{p+1})(\{p+1\}-\{m+n+s\})c_{m,n}=0.\label{odd-8}
    \end{align}

 We first compute $a_{m,n}$ with the assumption that $s\neq -1$. Taking $m=p=0$ and $n=s+1$ in \eqref{odd-2}, we obtain $\{s+1\}a_{0,s+1}=\{2s+2\}a_{0,s+1}$. Thus $a_{0,s+1}=0$. With this and setting $n=0$, $p=s+1$ in \eqref{odd-1}, we have $2\{m+2s+2\}a_{m,s+1}=(1+q^{s+1})\{m\}a_{m,s+1}$. It implies $a_{m,s+1}=0$ for all $m\in\Z$. Then we set $n=s+1$, $p=0$ in \eqref{odd-2} and obtain
    \begin{eqnarray*}
        (1+q^{s+1})(\{m+s+1\}-\{s+1\})a_{m,0}=0.
    \end{eqnarray*}
    It follows $a_{m,0}=0$ for $m\neq 0$. Taking $m=p=0$ in \eqref{odd-1} gives
    \[
        (1+q^n)(\{n\}-\{s+1\})a_{0,0}+2\{n+s+1\}a_{n,0}-2\{n\}a_{n,0}=0.
    \]
    We can always find an $n\in\Z$ such that $n\neq 0$ and $n\neq s+1$. Then we can deduce $a_{0,0}=0$ from the above equation. Thus $a_{m,0}=0,\ \  \forall\, m\in\Z.$ Finally setting $p=0$ in \eqref{odd-2} gives
       $ (1+q^m)\{n\}a_{m,n}=2\{m+n+s+1\}a_{m,n}.$
    It implies
    \begin{eqnarray*}
        a_{m,n}=0,\ \  \forall\ m,n\in\Z.
    \end{eqnarray*}

    When $s=-1$, Equations \eqref{odd-1} and \eqref{odd-2} become respectively
    \begin{eqnarray*}
        (1+q^n)(\{n\}-\{m+p\})a_{m,p}+(1+q^m)(\{n+p\}-\{m\})a_{n,p}-(1+q^{p})(\{n\}-\{m\})a_{m+n,p}=0,\\
        (1+q^n)(\{m+p\}-\{n\})a_{m,p}-(1+q^m)(\{p\}-\{n\})a_{m,n+p}+(1+q^{p})(\{p\}-\{m+n\})a_{m,n}=0.
    \end{eqnarray*}
    They are the same with \eqref{witt-1} and \eqref{witt-2} respectively, where $s=0$. Hence  $a_{m,n}=(\{n\}-\{m\})a_{0,1},~\forall~m,n\in \Z$.

    Next we compute $b_{m,n}=0$. Setting $n=p+1$ in \eqref{odd-4} gives
    \[
        (1+q^{p+1})\big(\{m+p+s\}-\{p+1\}\big) b_{m,p}=0.
    \]
    It implies $b_{m,p}=0$ for all $m\neq 1-s$ and $p\in\Z$. We can always find $m_0, n_0\in\Z$ such that $m_0\neq 0$, $n_0\neq 0$, $m_0\neq n_0$ and $m_0+n_0=1-s$. Setting $m=m_0$ and $n=n_0$ in \eqref{odd-3}, we have $(1+q^{p+1})(\{n_0\}-\{m_0\})b_{1-s,p}=0$. It follows that $b_{1-s,p}=0$ for all $p\in\Z$. Hence
      $b_{m,p}=0,~\forall~m,p\in\Z$.

    To compute $c_{m,n}=0$, we take $n=m+1$ in \eqref{odd-5} and obtain
    \[
        (1+q^{m+1})\big(\{m+1\}-\{m+p+s\}\big) c_{m,p}=0.
    \]
    This gives $c_{m,p}=0$ for $p\neq 1-s$ and $m\in\Z$. We can always find $p_0, n_0\in\Z$ such that $p_0\neq 0$, $n_0\neq 0$, $p_0\neq n_0$ and $p_0+n_0=1-s$. Then setting $p=p_0$ and $n=n_0$ in \eqref{odd-6},we get $c_{m,1-s}=0$, for all $m\in\Z$. Hence $c_{m,p}=0,~\forall~m,p\in\Z$.

Finally, by $c_{m,n}=0$ and using the same discussions to $b_{m,n}$, we can obtain $d_{m,n}=0,~\forall~m,n\in \Z$. This ends the proof.
\end{proof}

By Proposition \ref{prop6}, Theorem \ref{th3} (2) holds.

\subsection{Applications}

Now we apply the results obtained in previous subsections to determine linear commuting maps on $\mathcal {W}^q$, $\mathcal{V}^q$ and $\mathcal{L}^q$, respectively.

\begin{theo}\label{w22}
    Any linear commuting map $f$ on  $\mathcal {W}^q$ is exactly of the form
    \begin{align}\label{3.83}
        f(L_m)=\lambda L_m+\mu W_m,\,
        f(W_m)=\lambda W_m,~\forall~ m\in \Z.
    \end{align}
    for some $\lambda,~\mu \in \C$.
\end{theo}
\begin{proof}
    Let $f$ be a linear commuting map on $\mathcal {W}^q$. By Lemma \ref{lemma1}, we have a biderivation of $\mathcal {W}^q$ defined by
    \begin{align*}
        \varphi(x,y)=[f(x),y],~\forall~x,y \in \mathcal {W}^q.
    \end{align*}
    By Theorem \ref{th1},  $\varphi=\lambda\varphi_{ad}+\mu \varphi_0$ for some $\lambda, \mu \in \C$.
    More precisely,
    \begin{align*}
            \varphi(L_m,L_n) &=[f(L_m),L_n]=\lambda[L_m,L_n]+\mu[W_m,L_n], \\
            \varphi(L_m,W_n) &=[f(L_m),W_n]=\lambda[L_m,W_n],              \\
            \varphi(W_m,L_n) &=[f(W_m),L_n]=\lambda[W_m,L_n],              \\
            \varphi(W_m,W_n) &=[f(W_m),W_n]=0.
    \end{align*}
   Then \eqref{3.83} follows.
\end{proof}
\begin{coro}\label{coco1}
    \begin{itemize}
        \item [(1)] Let $f$  be a commuting automorphism of $\mathcal {W}^q$. Then $f$ is the identity map.
        \item[(2)] Let $f$  be a commuting derivation of $\mathcal {W}^q$. Then $f$ is the zero map.
    \end{itemize}
\end{coro}
\begin{proof}
    By Theorem \ref{w22}, we have $f(L_m)=\lambda L_m+\mu W_m,~f(W_m)=\lambda W_m$, for some $\lambda, \mu \in \C$.

    (1) If $f$  is a commuting automorphism of $\mathcal {W}^q$, then we have
    $f([L_m,L_n])=[f(L_m),f(L_n)].$
    It gives
    \begin{eqnarray*}
        \lambda [n-m]L_{m+n}+\mu [n-m]W_{m+n}=\lambda^2 [n-m]L_{m+n}+2\lambda \mu [n-m]W_{m+n}, \ \,\forall \ \, m,n\in\Z,
    \end{eqnarray*}
    and thus $\lambda=0,~\mu =0$ or $\lambda=1,~\mu =0$. It follows that $f$ is the zero map from the first case and $f$ is the identity map from the second case. Since $f$ is invertible, $f$ is the identity map.

    (2) If $f$  is a commuting derivation  of $\mathcal {W}^q$, then
    \begin{eqnarray*}
        f([L_m,L_n])=[f(L_m),L_n]+[L_m,f(L_n)], \ \,\forall \ \, m,n\in\Z.
    \end{eqnarray*}
    It follows
    \begin{eqnarray*}
        \lambda [n-m]L_{m+n}+\mu [n-m]W_{m+n}=2\lambda [n-m]L_{m+n}+2\mu [n-m]W_{m+n}, \ \,\forall \ \, m,n\in\Z.
    \end{eqnarray*}
    Hence $\lambda=0$ and $\mu=0$. It shows that $f$ is the zero map.
\end{proof}

Similarly, we have the following results for the $q$-deformed Witt algebra $\mathcal {V}^q$:
\begin{theo} Let $f$
    be a linear commuting map on $\mathcal {V}^q$. The following hold:
    \begin{itemize}
        \item [(1)]$f$ is a scalar
              multiplication map.
        \item [(2)] If $f$  be an automorphism of $\mathcal {V}^q$, then $f$ is the identity map.
        \item[(3)] If $f$  be a derivation of $\mathcal {V}^q$, then $f$ is the zero map.
    \end{itemize}
\end{theo}
\begin{proof}
    Similar to the proof of Theorem \ref{w22} and Corollary \ref{coco1}.
\end{proof}

For the $q$-deformed Witt superalgebra $\mathcal{L}^q$, we have the following:
\begin{theo}\label{thm:6.5}
    Any linear commuting map $f$ on $\mathcal{L}^q$ is exactly one of the following form:
    \begin{itemize}
        \item [(1)] If $f$ is even, then $f$ is a scalar
              multiplication map.
        \item [(2)] If $f$ is odd, then there exists some $\lambda \in \C$ such that
              \begin{align}
                  f(L_m)=\lambda G_{m-1},\
                  f(G_m)=0,~\forall~ m\in \Z.
              \end{align}
    \end{itemize}
\end{theo}
\begin{proof}
    Let  $f$ be a linear commuting map on $\mathcal{L}^q$. By Lemma \ref{lemma2},  we have a super-biderivation defined by
    \begin{align*}
        \varphi(x,y)=[f(x),y],~\forall~x,y \in \mathcal {L}^q.
    \end{align*}
    Obviously, $|f|=|\varphi|$.

    (1) Assume that $f$ is even. By Theorem \ref{th3} (1), we have $\varphi=\lambda\varphi_{ad}$ for some $\lambda\in \C$.
    More precisely,
    \begin{align*}
        \varphi(L_m,L_n) & =\lambda[L_m,L_n]=[\lambda L_m,L_n]=[f(L_m),L_n], \\
        \varphi(G_m,L_n) & =\lambda[G_m,L_n]=[\lambda G_m,L_n]=[f(G_m),L_n].
    \end{align*}
    It follows that $f(L_m)=\lambda L_m$ and $f(G_m)=\lambda G_m$ for all $m\in\Z$. Hence assertion (1) holds.

    (2) Assume that $f$ is odd. By Theorem \ref{th3} (2), we have $\varphi=\lambda\varphi_{-1}$ for some $\lambda\in \C$.
    Then we can easily deduce that $f(L_m)=\lambda G_{m-1}$ and $f(G_m)=0$ for all $m\in\Z$. Hence assertion (2) holds.
\end{proof}

\begin{coro}\begin{itemize}
        \item [(1)] Let $f$ be a commuting automorphism of $\mathcal {L}^q$. Then $f$ is the identity map.
        \item[(2)] Let $f$  be a commuting super-derivation of $\mathcal {L}^q$. Then $f$ is the zero map.
    \end{itemize}
\end{coro}
\begin{proof}
    (1) If $f$  is a commuting automorphism, then $f$ is an invertible even map. By Theorem \ref{thm:6.5}, $f$ is a scalar multiplication map, i.e., $f=\lambda {\rm id}$ for some $\lambda\in\C\setminus\{0\}$. Furthermore, it follows from
    \begin{eqnarray*}\label{eqn: 6.14}
        \lambda[L_m,L_n]=f([L_m,L_n])=[f(L_m),f(L_n)]=\lambda^2[L_m,L_n],
    \end{eqnarray*}
    that $\lambda=1$. Hence, $f$ is the identity map.

    (2) If $f$  is a commuting super-derivation, then
    \begin{eqnarray}\label{*}
        f([x,y])=[f(x),y]+(-1)^{|f||x|}[x,f(y)],\ \,\forall\ x\in\mathcal{H}(\mathcal {L}^q), \,y\in \mathcal {L}^q.
    \end{eqnarray}
    Taking the pair $(x,y)$ to be $(L_m,L_n)$ in \eqref{*}, we have
    \begin{eqnarray*}
        f([L_m,L_n])=[f(L_m),L_n]+[L_m,f(L_n)], \ \,\forall \ \, m,n\in\Z.
    \end{eqnarray*}
    By Theorem \ref{thm:6.5}, if $f$ is even, then
    \begin{eqnarray*}
        \lambda[L_m,L_n]=2\lambda[L_m,L_n], \ \,\forall \ \, m,n\in\Z.
    \end{eqnarray*}
    It follows that $\lambda=0$ and thus $f$ is the zero map. Whereas if $f$ is odd, we have
    \begin{align*}
        \lambda(\{n\}-\{m\})G_{m+n-1}=[\lambda G_{m-1},L_n]+[L_m,\lambda G_{n-1}]
        =2\lambda(\{n\}-\{m\})G_{m+n-1}, \ \,\forall \ \, m,n\in\Z.
    \end{align*}
    This gives $\lambda=0$. Hence $f$ is the zero map.

    The proof is complete.
\end{proof}


\section{$\alpha$-Biderivations}

In this section, we introduce the notions of $\alpha$-derivations and $\alpha$-biderivations for a Hom-Lie algebra $(\G,[\cdot,\cdot],\alpha)$, and we establish a close relation between $\alpha$-derivations and $\alpha$-biderivations. The superanalogue is also discussed. As an application, we will show that the $q$-deformed $W(2,2)$-algebra, the $q$-deformed Witt algebra and superalgebra have no nonzeo $\alpha$-biderivations or $\alpha$-super-biderivations.

Let $(\G,[\cdot,\cdot],\alpha)$ be a Hom-Lie algebra or a Hom-Lie superalgebra. For any nonnegative integer $k$, denote by $\a^k$
the $k$-times composition of $\a$, i.e. $\a^k = \a\circ \cdots \circ \a$ ($k$-times).  We introduce the following concepts.

\begin{defi}\label{def1}{\rm
        \begin{itemize}
            \item[(1)]
        A linear map $D: \mathscr{G} \to \G$  is called an $\alpha^k$-derivation of a Hom-Lie algebra $(\G,[\cdot,\cdot],\alpha)$,  if it satisfies
        \begin{eqnarray}\label{7-1}
            D([x,y])=[D(x),\alpha^k(y)]+[\alpha^k(x),D(y)],~\forall~x,y\in \G.
        \end{eqnarray}
        \item[(2)]  A homogeneous linear map $D:\mathscr{G}\rightarrow \mathscr{G}$ is called an $\alpha^k$-super-derivation of a Hom-Lie superalgebra $(\mathscr{G},[\cdot,\cdot],\alpha)$  if it satisfies
    $$D([x,y])=[D(x),\alpha^k(y)]+(-1)^{|x||D|}[\alpha^k(x),D(y)], \ \,\forall\ \, x\in \mathcal{H}(\mathscr{G}),\ \, y\in \mathscr{G}. $$
      \end{itemize}
    }
\end{defi}

\begin{remark}\label{re1} The notion of an $\alpha^k$-derivation of Hom-Lie algebras was given in \cite[Definition 3.1]{SH}, and the superanalogue was given in \cite[Definition 2.1]{AMS1}. Compared with them, the strong assumption $D\circ \alpha=\alpha\circ D$ is removed in Definition \ref{def1}.
\end{remark}

\begin{defi}{\rm
        \begin{itemize}
            \item[(1)] An $\a$-biderivation of a Hom-Lie algebra $(\mathscr{G},[\cdot,\cdot],\alpha)$ is a bilinear map $\varphi:\mathscr{G}\times \mathscr{G}\rightarrow \mathscr{G}$ satisfying
                \begin{align}
                    \varphi([x,y],z)= & ~[\varphi(x,z),\alpha(y)]+[\alpha(x),\varphi(y,z)],\label{vf1-1} \\
                    \varphi(x,[y,z])= & ~[\varphi(x,y),\alpha(z)]+[\alpha(y),\varphi(x,z)],\label{vf2-1}
                \end{align}
                for all $x,y,z\in \mathscr{G}$.
            \item[(2)]  An $\a$-super-biderivation of a Hom-Lie superalgebra $(\mathscr{G},[\cdot,\cdot],\alpha)$ is a homogeneous bilinear map $\varphi:\mathscr{G}\times \mathscr{G}\rightarrow \mathscr{G}$ satisfying \begin{align}
                    \varphi \big([x,y],z\big)=~ & (-1)^{|y||z|}[\varphi(x,z),\alpha(y)]+(-1)^{|\varphi||x|}[\alpha(x),\varphi(y,z)],\label{svf1-1} \\
                    \varphi\big(x,[y,z]\big)=~  & [\varphi(x,y),\alpha(z)]+(-1)^{(|\varphi|+|x|)|y|}[\alpha(y),\varphi(x,z)],\label{svf2-2}
                \end{align}
                for all homogeneous elements $x, y, z$ in $\mathscr{G}$.
        \end{itemize} }
\end{defi}

In the following lemma, we characterize a basic relationship between $\a$-derivations and $\a$-biderivations.
\begin{lemm}\label{lemma3} Let $\varphi$ be an $\a$-biderivation of a Hom-Lie algebra $(\mathscr{G},[\cdot,\cdot],\alpha)$. For every $z\in\mathscr{G}$,
    the maps $x\mapsto\varphi(x,z)$ and  $x\mapsto \varphi(z,x)$ are $\a$-derivations of\, $\mathscr{G}$.
\end{lemm}
\begin{proof} For convenience, we write the map $x\mapsto\varphi(x,z)$ as
    $f(x)=\varphi(x,z),~\forall x \in \G$. We only need to verify that $f$ satisfies \eqref{7-1} with $k=1.$ For all $x,y\in \mathscr{G}$, by \eqref{vf1-1} we have
    \begin{align*}
        f([x,y])=\varphi([x,y],z)=[\varphi (x,z),\alpha(y)]+[\alpha(x), \varphi(y,z)]
        =[f(x),\alpha(y)]+[\alpha(x),f(y)].
    \end{align*}
    It follows from \eqref{7-1} that $f$ is an  $\alpha$-derivation. Similarly, we can prove that $\varphi$ is also
    an $\a$-derivation with respect to the second component.
\end{proof}

The supercase of Lemma \ref{lemma3} is also valid with a similar proof.
\begin{lemm}\label{lemma4} Let $\varphi$ be an $\a$-super-biderivation of a Hom-Lie superalgebra $(\mathscr{G},[\cdot,\cdot],\alpha)$.  For every $z\in\mathscr{G}_{\bar 0}$,
    the maps $x\mapsto\varphi(x,z)$ and  $x\mapsto \varphi(z,x)$ are $\a$-super-derivations of\, $\mathscr{G}$.
\end{lemm}

The following lemma is straightforward.
\begin{lemm}\label{lemma5}
    \begin{itemize}
        \item[(1)] Suppose that $(\mathscr{G},[\cdot,\cdot],\alpha)$ is Hom-Lie algebra without nonzero $\a$-derivations. Then $\mathscr{G}$ has no nonzero $\a$-biderivations.
        \item[(2)] Suppose that $(\mathscr{G},[\cdot,\cdot],\alpha)$ is Hom-Lie superalgebra without nonzero $\a$-super-derivations and
            $\varphi$ is an $\a$-super-biderivation of $\mathscr{G}$. Then, for every $z\in\mathscr{G}_{\bar 0}$, we have
            \begin{eqnarray*}
                \varphi(z,x)=\varphi(x,z)=0,\ \forall\ \, x\in \mathcal{H}(\mathscr{G}).
            \end{eqnarray*}
    \end{itemize}
\end{lemm}

As we have pointed out in Remark \ref{re1}, the difference between
the definitions of $\a$-(super-)derivations given in \cite{SH,AMS1} and our definition is
the assumption $D\circ \alpha=\alpha\circ D$. However, we find that the authors in \cite{AMS1,YY} did not use the condition $D\circ \alpha=\alpha\circ D$ when computing $\a$-(super-)derivations of the $q$-deformed $W(2,2)$ algebra and Witt superalgebra. Although the condition $D\circ \alpha=\alpha\circ D$ could make their calculations much simpler, it had no impact on  the final results. All in all, the following results due to \cite{AMS1,YY} are also right in our framework of this paper.

\begin{prop}\label{prop3} \begin{itemize} \item[(1)] Both $\mathcal{W}^q$ and $\mathcal{V}^q$ have no nonzero $\a$-derivations.
    \item[(2)] $\mathcal{L}^q$ has no nonzero $\a$-super-derivations.
   \end{itemize}
\end{prop}

Combining Lemma \ref{lemma5} with Proposition \ref{prop3}, we have the following.

\begin{prop}\label{prop4} \begin{itemize}
    \item[(1)] Both $\mathcal{W}^q$ and $\mathcal{V}^q$ have no nonzero $\a$-biderivations.
    \item[(2)] $\mathcal{L}^q$ has no nonzero $\a$-super-biderivations.
   \end{itemize}
\end{prop}
\begin{proof} (1) It follows immediately from Lemma \ref{lemma5} (1) and Proposition \ref{prop3} (1).

(2) Let $\varphi$ be an $\a$-super-biderivation on $\mathscr{L}^q$ of degree $s$. By Lemma \ref{lemma5} (2), we have
    \begin{eqnarray}\label{**}
        \varphi(L_m,L_n)=\varphi(L_m,G_n)=\varphi(G_m,L_n)=0, \, \forall \, m,n\in\Z.
    \end{eqnarray}
    It remains to evaluate the values of $\varphi$ on the pairs of $(G_m,G_n)$ for all $m,n\in\Z.$ Taking the triple $(x,y,z)$ to be $(G_m,G_n,L_p)$ in \eqref{svf2-2} and by \eqref{**}, we obtain
    \begin{eqnarray}\label{***1}
        \varphi(G_m,[G_n,L_p])=[\varphi(G_m,G_n),\a(L_p)].
    \end{eqnarray}
    Taking the triple $(x,y,z)$ to be $(G_m,G_n,G_p)$ and $(L_n,G_m,G_p)$ respectively in \eqref{svf1-1} and by \eqref{**}, we obtain
    \begin{align}
        [\varphi(G_m,G_p),\a(G_n)]= & ~(-1)^{|\varphi||G_m|}[\a(G_m),\varphi(G_n,G_p)],\label{***2} \\
        \varphi([L_n,G_m],G_p)=     & ~[\a(L_n),\varphi(G_m,G_p)].\label{***3}
    \end{align}

    If $|\varphi|=\bar{0}$, we can assume that $\varphi(G_m,G_n)=a_{m,n}L_{m+n+s}$, for all $m,n\in\Z.$  Inserting this into \eqref{***1} and \eqref{***2}, respectively, gives
    \begin{align}
        (\{p\}-\{n+1\})a_{m,n+p}=              & ~(1+q^p)(\{p\}-\{m+n+s\})a_{m,n},\label{**2}       \\
        (1+q^{n+1})(\{n+1\}-\{m+p+s\})a_{m,p}= & ~(1+q^{m+1})(\{n+p+s\}-\{m+1\})a_{n,p}\label{**3}.
    \end{align}
    Setting $n=0$ and $p=1$ in \eqref{**2} gives $(1-\{m+s\})a_{m,0}=0$, and thus $a_{m,0}=0$ for $m\neq 1-s$.
    Taking $p=0$ and $m=1-s$ in \eqref{**3} gives
    \begin{eqnarray*}
        (1+q^{n+1})(\{n+1\}-1)a_{1-s,0}=(1+q^{2-s})(\{n+s\}-\{2-s\})a_{n,0}.
    \end{eqnarray*}
    In this equation, we can always find a suitable $n$ such that $n\neq 1-s$ and $n\neq 0$ and hence $a_{1-s,0}=0$. Now we have $a_{m,0}=0$ for all $m\in\Z$. Taking this into account and setting $n=0$ in \eqref{**2}, we obtain $(\{p\}-1)a_{m,p}=0$ and thus $a_{m,p}=0$ for all $m\in\Z$ and $p\in\Z\setminus\{1\}$. With this and letting $p=2$, $n=-1$ in \eqref{**2} we have $a_{m,1}=0$. Hence, $a_{m,n}=0$ for all $m,n\in\Z$.

    If $|\varphi|=\bar{1}$, we can assume that $\varphi(G_m,G_n)=b_{m,n}G_{m+n+s}$, for all $m,n\in\Z.$ Inserting this into \eqref{***1} and \eqref{***3} gives
    \begin{align}
        (\{p\}-\{n+1\})b_{m,n+p}= & ~(1+q^p)(\{p\}-\{m+n+s+1\})b_{m,n},\label{**4} \\
        (\{m+1\}-\{n\})b_{m+n,p}= & ~(1+q^n)(\{m+p+s+1\}-\{n\})b_{m,p}.\label{**5}
    \end{align}
    By exchanging the index $p$ and $n$ in \eqref{**5}, we obtain
    \begin{eqnarray}
        (\{m+1\}-\{p\})b_{m+p,n}=(1+q^p)(\{m+n+s+1\}-\{p\})b_{m,n}.\label{**6}
    \end{eqnarray}
    Comparing \eqref{**4} and \eqref{**6}, we have
    \begin{eqnarray}
        (\{m+1\}-\{p\})b_{m+p,n}=(\{n+1\}-\{p\})b_{m,n+p}.\label{**7}
    \end{eqnarray}
    Taking $p=0$ in \eqref{**7} gives $b_{m,n}=0$ for $m\neq n$. By using this and setting $p=1$, $n=m-1$ in \eqref{**7}, we have $b_{m,m}=0$ for $m\neq 1$. Finally, taking $m=1,n=2,p=-1$ in \eqref{**7}, we obtain $b_{1,1}=0$. This proves $b_{m,n}=0$ for all $m,n\in\Z$.

    Combining the two cases above, we get the result.
\end{proof}

At the end of this section, we present an example of Hom-Lie superalgebras with nonzero $\a$-super-derivations and $\a$-super-biderivations.

\begin{ex}{\rm Let $L=L_{\bar 0}\bigoplus L_{\bar 1}$ be a three-dimensional superspace equipped with the following nonzero products
        $$[x_1, x_2] = \lambda^2 x_1,\ [x_2, y] =- \frac{\lambda}{2} y,\ [y, y] =\lambda^2 x_1,$$
        where $L_{\bar0} = \C x_1\bigoplus \C x_2, L_{\bar1} = \C y$, and $\lambda\in\C\setminus\{0,1\}$.  Define a linear even map $\a: L\rightarrow L$ by
        $$\a (x_1) = \lambda^2x_1,\ \a(x_2) = x_2,\ \a(y) = \lambda y.$$
        Then $(L,[\cdot,\cdot],\a)$ is a Hom-Lie superalgebra  (see \cite{ZHB}). Define a linear map $D:L\rightarrow L$ by
        \begin{eqnarray*}
            D(x_1)=2c\lambda x_1,\ D(x_2)=bx_1,\ D(y)=cy,\ \ b,c\in\C.
        \end{eqnarray*}
         It is easy to check that $D$ is a nonzero even $\a$-super-derivation of $L$ if and only if $(b,c)\neq (0,0)$.  Moreover, define a bilinear map $\varphi:L\times L\rightarrow L$ as follows:
        \begin{alignat*}{3}
             & \varphi(x_1,x_1)=0,\quad &  & \varphi(x_1,x_2)=-2d\lambda x_1,~~~&  & \varphi(x_1,y)=0,   \\
             & \varphi(x_2,x_2)=a x_1,\quad &  & \varphi(x_2,x_1)=2d\lambda x_1,~~  &  & \varphi(x_2,y)=d y, \\
             & \varphi(y,x_1)=0,\quad   &  & \varphi(y,x_2)=-dy,~~   &  & \varphi(y,y)=kx_1,
        \end{alignat*}
        where $d=\frac{1}{\lambda^2}(\frac12-\lambda)k$, and $a,k\in\C$ . One can check that $\varphi$ is a nonzero even $\a$-super-biderivation of $L$ if and only if $(a,k)\neq (0,0)$, and $\varphi$ is not skew-symmetric when $a\neq 0$.
    }\end{ex}

 \vs{5pt}
\noindent{\bf{ Acknowledgements.}}\ {Supported by National Natural Science
Foundation grants of China (11301109).



\end{document}